\documentclass[a4paper,11pt]{article}
\usepackage{amsmath}
\usepackage{amsthm}
\usepackage{amssymb}
\usepackage{amsfonts}
\usepackage{graphicx}
\usepackage{enumerate}
\usepackage[all]{xy}
\usepackage{epsfig}

\normalsize

\newdir{ >}{!/-8.5pt/@{>}}

\newtheoremstyle{theoremstyle}
  {10pt}      
  {5pt}       
  {\itshape}  
  {}          
  {\bfseries} 
  {:}         
  {.5em}      
  {}          

\newtheoremstyle{examplestyle}
  {10pt}      
  {5pt}       
  {}          
  {}          
  {\bfseries} 
  {:}         
  {.5em}      
  {}          

\theoremstyle{theoremstyle}
\newtheorem{theorem}{Theorem}[section]
\newtheorem*{theorem*}{Theorem}
\newtheorem{lemma}[theorem]{Lemma}
\newtheorem{proposition}[theorem]{Proposition}
\newtheorem*{proposition*}{Proposition}
\newtheorem{corollary}[theorem]{Corollary}
\newtheorem*{corollary*}{Corollary}

\theoremstyle{examplestyle}
\newtheorem{example}[theorem]{Example}
\newtheorem{definition}[theorem]{Definition}
\newtheorem{definition*}{Definition}
\newtheorem{remark}[theorem]{Remark}
\newtheorem{remark*}{Remark}

\newcommand{\comment}[1]{}

\newcommand{\rays}{{\Delta(1)}}

\newcommand{\weildivisors}{{\mathbb{Z}^\rays}}
\newcommand{\mweildivisors}{$\weildivisors$}
\newcommand{\sh}[1]{\mathcal{#1}}
\newcommand{\msh}[1]{$\sh{#1}$}
\newcommand{\spec}[1]{\operatorname{spec}(#1)}

\newcommand{\orb}[1]{\operatorname{orb}(#1)}

\newcommand{\codim}{\operatorname{codim}}

\newcommand{\rk}{\operatorname{rk}}
\newcommand{\ob}{\operatorname{Ob}}
\newcommand{\mor}{\operatorname{Mor}}
\newcommand{\Hom}{\operatorname{Hom}}

\newcommand{\kvect}{k\operatorname{-\bf Vect}}

\newcommand{\lcm}{\operatorname{lcm}}

\newcommand{\zip}{\operatorname{zip}}
\newcommand{\unzip}{\operatorname{unzip}}

\newcommand{\Z}{\mathbb{Z}}
\newcommand{\R}{\mathbb{R}}
\newcommand{\N}{\mathbb{N}}
\newcommand{\p}{\mathcal{P}}
\newcommand{\cL}{\mathcal{L}}
\newcommand{\n}{{\underline{n}}}

\title{Resolutions for Equivariant Sheaves over Toric Varieties}

\author{Markus Perling\footnote{Institut f\"ur Mathematik, Universit\"at
Paderborn, 33098 Paderborn, Germany, {\tt perling@math.upb.de}}}

\date{July 2004}

\begin{document}

\maketitle

\begin{abstract}
In this work we construct global resolutions for general coherent equivariant sheaves
over toric varieties. For this, we use the framework of sheaves
over posets. We develop a notion of gluing of posets and of sheaves over posets,
which we apply to construct global resolutions for equivariant sheaves. Our
constructions
give a natural correspondence between resolutions for reflexive equivariant sheaves
and  free resolutions of vector space arrangements.
\end{abstract}

\newpage

\tableofcontents

\newpage

\section{Introduction}

An important part of the theory of vector bundles over homogeneous spaces $G/P$ is
the study of {\em homogeneous} vector bundles. This class of vector bundles has first
been investigated by Kostant and Bott in the 50's, who clarified the relation
between the representation theory of $G$ and homogeneous bundles. This relation
in many cases allows to determine properties of homogeneous vector bundles very
explicitly, and so homogeneous bundles have played a great role in the field of
studying general vector bundles, notably over projective spaces.

In a more general situation, one considers a {\em quasi-homogeneous} space, i.e. a
space $X$ together with the action of an algebraic group $G$ such that this action
has a dense open orbit in $X$. In this context it is customary to speak about {\em
equivariant} rather than homogeneous vector bundles; denote $\sigma, p_2 : G \times
X \longrightarrow X$ the group action and the projection onto the second factor,
respectively, then a vector bundle (or a more general sheaf) \msh{E} on $X$ is {\em
equivariant} if there exists an isomorphism
\begin{equation*}
\Phi: \sigma^* \sh{E} \overset{\cong}{\longrightarrow} p_2^* \sh{E}
\end{equation*}
such that
\begin{equation*}
(\mu \times 1_X)^* \Phi = p_{23}^*\phi \circ (1_G \times \sigma)^* \Phi,
\end{equation*}
where $\mu$ is the group multiplication morphism and $p_{23}$ the projection onto the
second and third factor of $G \times G \times X$ (see also \cite{GIT}).
This situation in general is considerably more difficult than the case of homogeneous
spaces, as (at least) the following two things can happen: in general, $X$ has a
rather complicated
orbit structure, such that there are lower-dimensional invariant loci which allow
equivariant
vector bundles to degenerate to more general equivariant sheaves, if considered in
families in a suitable sense; moreover, the representation theory of $G$ contributes
only marginal information. So the conclusion is that one has to study the
complete category of equivariant sheaves over $X$, which in particular means:
\begin{enumerate}[(i)]
\item construct good invariants for equivariant sheaves over $X$,
\item study moduli spaces with respect to these invariants.
\end{enumerate}
In this work, we attempt to carry out part of such a program for equivariant sheaves
over toric varieties, which are probably the easiest examples of quasi-homogeneous
spaces.

\paragraph{Reflexive Sheaves.}
 Our approach is based on the framework of $\Delta$-families which we have
developed in earlier work (\cite{perling1}), which in turn
generalizes the characterization of Klyachko (\cite{Kly90}, \cite{Kly91}) of
equivariant reflexive sheaves. Klyachko's observation was that every such sheaf
\msh{E} is equivalent to a finite dimensional vector space $\mathbf{E}$ together
with a finite set of full filtrations
\begin{equation*}
\cdots \subset E^\rho(i) \subset E^\rho(i + 1) \subset \cdots \subset \mathbf{E}
\end{equation*}
for $i \in \Z$ and every torus invariant divisor $\rho \in \rays$ (see section
\ref{toricvarieties} for notation). Naively, one can seperate two kinds of data from
such a set of filtrations: first, the indices $i$, preferably those where the
dimension of the filtration jumps $E^\rho(i) \subsetneq E^\rho(i + 1)$, and second,
the flags underlying the filtrations, when we forget about the indices. One could
think of the indices as a discrete invariant for \msh{E}, and the flags as {\em
moduli} for the sheaf. However, it turns out that the indices essentially only
determine the first equivariant Chern class of \msh{E}, and the moduli of flags do
not behave very well in sheaf theoretic sense. This has been investigated in detail
in \cite{perling2} for case of equivariant vector bundles of rank two over toric
surfaces.

One could proceed now and declare the equivariant Chern classes as invariants for
equivariant sheaves and construct moduli with respect to these (this has been done in
\cite{perling2}), but we are interested in a more direct approach and want to analyze
the flags underlying the filtrations. These flags and their intersections determine
a subvector space arrangement of $\mathbf{E}$, and as there is no more data left to
describe \msh{E}, one intuitively assumes that all further properties of \msh{E} are
somehow encoded in this arrangement.

Our approach is to construct a global resolution
for any given equivariant sheaf \msh{E} over $X$. From the point of view of
homogeneous coordinate rings (see \cite{Cox}) it has been observed (\cite{ERR}) that
every such sheaf has a finite global resolution
\begin{equation*}
\label{resolutionsequence}
0 \longrightarrow \sh{F}_s \longrightarrow \cdots \longrightarrow \sh{F}_0
\longrightarrow \sh{E} \longrightarrow 0
\end{equation*}
where $\sh{F}_i \cong \bigoplus_j \sh{O}(D_{ij})$ for every $i$. Here, the $D_{ij}$
are torus invariant Weil divisors, and the sheaves $\sh{O}(D_{ij})$ are equivariant
reflexive sheaves of rank one; in the case where $X$ is smooth, these sheaves always
are invertible.
We will give an explicit construction for such resolutions, which for the case of
reflexive sheaves will only depend on the underlying vector space arrangements.
Our results generalize a result of Klyachko, who in \cite{Kly90} constructed a
{\em canonical} resolution in the case where \msh{E} is locally free and $X$ is
smooth and complete.

\paragraph{Vector space arrangements.}
An interesting aspect of our construction is the solution of the following problem;
consider any subvector space arrangement in some vector space $\mathbf{E}$, and its
underlying poset $\p$ which is given by the set of subvector spaces in the
arrangement together with the partial order which is given by inclusion. Then, does
there exist a vector space $\mathbf{F}$ together with a {\em coordinate} vector space
arrangement such that the underlying poset is isomorphic to $\p$? The answer is yes,
and it is rather straightforward to see that one just needs to choose $\mathbf{F}$
large enough, such that the combinatorics of $\p$ can be modelled by coordinate
spaces of $\mathbf{F}$. As a byproduct, we obtain a surjection $\mathbf{F}
\twoheadrightarrow \mathbf{E}$ such that for every element $V \in \p$ and its
corresponding subvector space $F_V$ of $\mathbf{F}$, we have a commutative
exact diagram
\begin{equation*}
\xymatrix{
0 \ar[r] & \mathbf{K} \ar[r] & \mathbf{F} \ar[r] & \mathbf{E} \ar[r] & 0 \\
0 \ar[r] & K_V \ar[r] \ar@{^{(}->}[u] & F_V \ar[r] \ar@{^{(}->}[u] & V \ar[r]
\ar@{^{(}->}[u] & 0.
}
\end{equation*}
The vector spaces $K_V$ again form a vector space arrangement in $\mathbf{K}$ whose
underlying poset is a subset of the original poset $\p$. We call the arrangement
$K_V$ the {\em first syzygy arrangement} of $\p$.

By iterating this procedure, we obtain an exact sequence of vector spaces
\begin{equation*}
0 \longrightarrow \mathbf{F}_s \longrightarrow \cdots \longrightarrow \mathbf{F}_0
\longrightarrow \mathbf{E} \longrightarrow 0
\end{equation*}
where every $\mathbf{F}_i$ contains a coordinate vector space arrangement whose
underlying poset coincides with the poset underlying the $i$-th
syzygy arrangement. In the case where the arrangement in $\mathbf{E}$ is closed
under performing intersections, we have even a good notion of {\em minimal}
resolutions; we obtain a unique representation of such an arrangement in terms of
the purely combinatorial information encoded in the successive coordinate space
arrangements. In a sense, we can think of the resolution as providing a
``K-theory''-class in a suitable category of vector space arrangements. We formulate
the following

\

\noindent
{\bf Conjecture:} Let \msh{E} be a reflexive equivariant sheaf over a toric variety
$X$, then every property of \msh{E} depends only on the indices of the filtrations
$E^\rho(i)$ and the $K$-theory class of the underlying vector space arrangement.

\

\noindent
One can read this conjecture also the way that the ``K-theory''-class of a
vector space arrangement is its finest possible invariant. The class of coordinate
vector space arrangements is a well-studied subject (see \cite{BuchstaberPanov1}),
and it would be interesting to see whether properties of general arrangements can
be studied through free resolutions.

\paragraph{Poset representations.}
The construction of some global resolution for an arbitrary equivariant coherent sheaf
over $X$ is not necessarily a difficult task, but in general the organization of all
the needed data is rather elaborate. Any nuts and bolts approach, starting from
scratch, would probably be rather cumbersome for the reader to follow; therefore we
adopt in
this paper a more formal approach, by developing a certain amount of
framework in the context of poset representations. Such representations, as a subtopic
of quiver representations \cite{Gabriel72}, have been studied since long (see
\cite{Nazarova80}). Any poset $\p$ with a partial order $\leq$ in a natural way is
equivalent to some category. In this category the objects are the elements of $\p$,
and the morphisms are the relations $x \leq y$, i.e. there exists at most one
morphism between two objects $x, y \in \p$. A {\em representation} of $\p$ is a
functor $F: \p \longrightarrow \kvect$, $x \mapsto F_x$, the category
of vector spaces over some field $k$. The representations
themselves form an abelian category whose morphisms are the natural transformations.

On a poset $\p$ there exists a natural topology, which is generated by the basis
$U(x) = \{x \leq y \in \p\}$. Using this topology, every representation $F$ of $\p$
induces a sheaf over $\p$ by setting $\sh{F}\big(U(x)\big) := F_x$, and conversely,
any sheaf over $\p$ with values in $\kvect$ induces a representation of $\p$. In fact,
the categories of representations of $\p$ and of sheaves over $\p$ with values in
$\kvect$ are equivalent. However, it will be more comfortable for us to have both
points of view in mind and to switch the picture freely. The additional bonus of
sheaves over $\p$ is that by the continuation to the whole topology of $\p$, they
automatically incorporate inverse limits via $\sh{F}(U) =
\underset{\leftarrow}{\lim} \sh{F}\big(U(x)\big)$, where the limit runs over all
$x \in U$. For
us, this is a very natural way to encode all possible pullback diagrams over the poset
$\p$. Sheaves over posets have been in the literature before, the first reference
we are aware of being \cite{Baclawski}. More recently, this kind of sheaves has been
used in similar contexts like ours, for the study of certain modules over semigroup
rings \cite{Yanagawa01}, and for vector space arrangements
\cite{DeligneGoreskyMacPherson}.

\paragraph{Posets and graded modules.}
Our general principle will be to start with local constructions and to globalize these
by some gluing procedure, where 'local' and 'global' means over affine and general
toric varieties, respectively. Recall that an affine toric variety $U_\sigma$ over
some algebraically
closed field $k$ on which the torus $T$ acts, is equivalent to the spectrum of a
normal semigroup ring $k[\sigma_M]$. The semigroup $\sigma_M$ is a subsemigroup of 
the character group $M \cong \Z^r$ of $T$, which is given by the intersection of a
convex rational polyhedral cone $\check{\sigma}$ in $M \otimes_\Z \R$ with $M$.
Any equivariant sheaf \msh{E} over an affine toric variety $U_\sigma$ is equivalent
to an $M$-graded $k[\sigma_M]$-module $E^\sigma = \Gamma(U_\sigma, \sh{E})$, i.e.
\begin{equation*}
E^\sigma = \bigoplus_{m \in M} E^\sigma_m.
\end{equation*}
A fundamental observation is that this grading is the reason that equivariant sheaves
over toric varieties have still a {\em semi-combinatorial} nature, in contrast to the
completely combinatorial description of the toric varieties themselves. To see this,
note that $\sigma_M$ endows $M$ with the structure of a poset by setting $m
\leq_\sigma m'$ iff $m' - m \in \sigma_m$ (we simplify here, as in fact this in
general only defines a preorder). This way, $E^\sigma$ is equivalent to a
representation of $M$ which maps every $m$ to the vector space $E^\sigma_m$, and every
relation $m \leq_\sigma m'$ is mapped to the vector space homomorphism $E^\sigma_m
\longrightarrow E^\sigma_{m'}$, which is given by multiplication with the monomial
$\chi(m' - m)$. It turns out that the category of representations of the {\em poset}
$M$ is equivalent to the category of equivariant quasicoherent sheaves over
$U_\sigma$.

\begin{figure}[htb]
\begin{center}
\includegraphics[height=4cm]{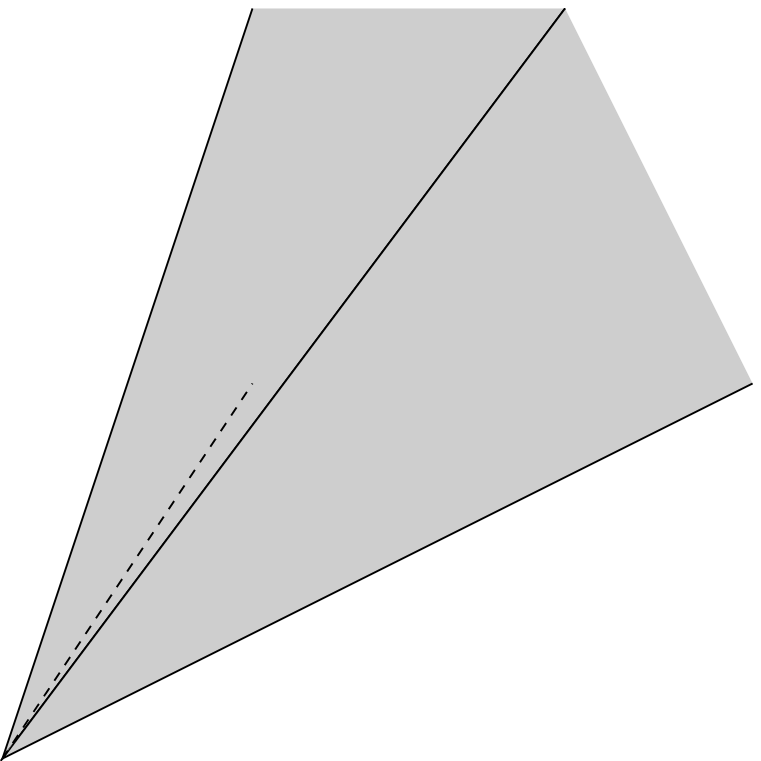}\qquad\qquad
\includegraphics[height=4cm]{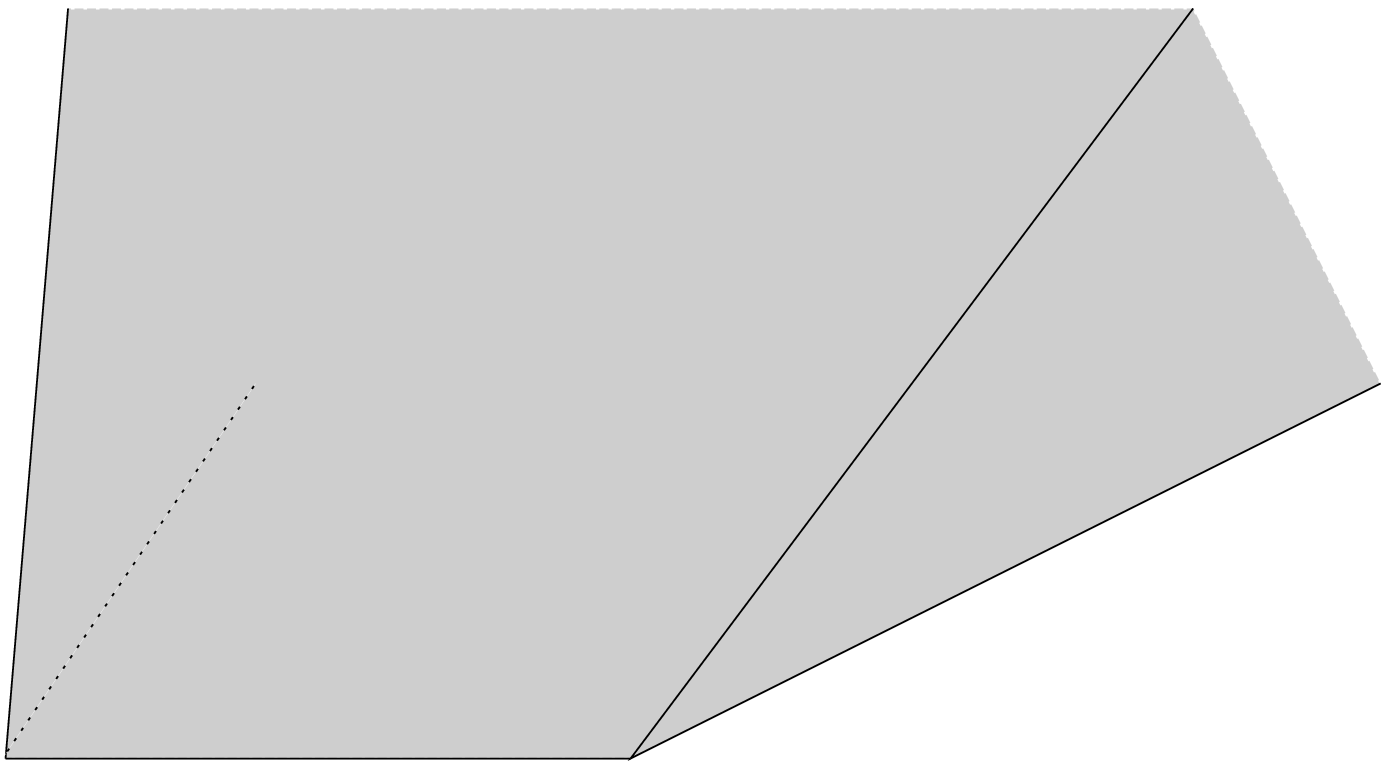}
\end{center}
\caption{Example of a cone with four maximal faces and a polyhedron}\label{f-introexample2}
\end{figure}

\begin{figure}[htb]
\begin{center}
\includegraphics[height=4cm,width=4cm]{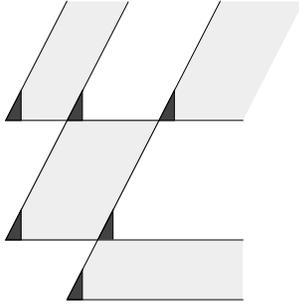}
\end{center}
\caption{A polyhedral decomposition of $\Z^2$ into seven regions}\label{f-introexample1}
\end{figure}

To be able to work truly with a finitely generated module, one needs an expedient
finite representation for it. For this, we introduce the notion of a {\em polyhedral
decomposition} of $M$. For any $\rho \in \sigma(1)$ and any integer $n_\rho$ we have
the shifted halfspace $\{m \in M_\R \mid \langle m, n(\rho) \rangle \geq n_\rho\}$ in
$M_\R$, and
for any tuple $\n = \big(n_\rho \mid \rho \in \sigma(1)\big) \in \Z^{\sigma(1)}$ the
intersection of half spaces $P_\n = \{m \mid \langle m, n(\rho) \rangle \geq n_\rho
\text{ for all } \rho \in \sigma(1)\}$. We call such an unbounded domain $P_\n$ a
{\em polyhedron}. Note that the dual cone $\hat{\sigma}$ itself is a polyhedron which
has the zero face as its unique compact face. Figure \ref{f-introexample2} shows an
example of a cone where $\sigma(1)$ consists of four rays, and a polyhedron defined
with respect to these four rays.
The intersection of two polyhedra $P_{\n_1}, P_{\n_2}$ is again a polyhedron, $P_\n$,
where $\n = (\max\{(n_{1, \rho}, n_{2, \rho}\} \mid \rho \in \sigma(1))$. This way,
any collection of polyhedra $P_{\n_1}, \dots, P_{\n_s}$ gives rise to a partition
of $M$ as follows. Define the 'least common multiple' $\n$ of any collection
$\n_{i_1}, \dots, \n_{i_k}$, by the componentwise maximum of the $\n_{i_j}$. Then the
equivalence classes $T_\n$ contain all $m \in M$ with $\langle m, n(\rho) \rangle \geq
n_{i, \rho}$ for all $\rho \in \sigma(1)$, for which there is no bigger least common
multiple $\n'$ satisfying these inequalities. Figure \ref{f-introexample1} shows a
partition of $\Z^2$ generated by three polyhedra.
The set of $\lcm$'s of the $\n_1, \dots, \n_s$ in a natural way becomes a poset, as a
subposet of $\Z^{\sigma(1}$ with partial order induced by the componentwise order. The
$\lcm$'s are a special case of a polyhedral decomposition which is induced by an
{\em admissible poset}. A finite poset $\p^\sigma \subset \Z^{\sigma(1)}$ is
admissible if for any $m \in M$ there exists a unique maximal element $\n \in
\p^\sigma$ such that $\langle m, n(\rho) \rangle \geq n_\rho$ for all $\rho \in
\sigma(1)$. $\p^\sigma$ is admissible {\em with respect to $E^\sigma$} if moreover
for every $\n \in \p^\sigma$ there exist a vector space $E_\n$ such that
$E_\n \cong E_m^\sigma$ for all $m \in T_\n$. The vector spaces $E_\n$ together with
appropriate morphisms $E_\n \longrightarrow E_{\n'}$ (whose existence is part of our
definition \ref{admissibledef} for admissible posets),
yield a representation of $\p^\sigma$, which encodes the
complete structure of $E^\sigma$. We can think of it, euphemistically, as a
compression of $E^\sigma$.

The most important feature of our constructions is that the compression of $E^\sigma$
is functorial, because we systematically exploit the
formalism of sheaves on posets; we finally arrive at an equivalence of categories
between sheaves over $\p^\sigma$ and $k[\sigma_M]$-modules with respect to which
$\p^\sigma$ is admissible (theorem \ref{admissibleequivalence}).
This in particular enables us to construct resolutions of $E^\sigma$ in terms of
free resolutions of the sheaf $E_\n$ over $\p^\sigma$. The resolutions obtained this
way are not free resolutions, but rather resolutions by reflexive modules of rank one.
Any $\n \in \Z^{\sigma(1)}$ gives rise to a $T$-invariant Weil divisor $D_\n = -
\sum_{\rho \in \sigma(1)} n_\rho D_\rho$ on $U_\sigma$, and thus to a reflexive sheaf
of rank one $\sh{O}_{U_\sigma}(D_\n)$. Write $S_{(\n)}$ for the associated reflexive
$k[\sigma_M]$-module, then its $M$-graded decomposition is given by
\begin{equation*}
S_{(n)} \cong \bigoplus_{m \in P_\n \cap M} k \cdot \chi(m).
\end{equation*}
Every equivalence class $T_\n$ has the shape of the forepart of the polyhedron $P_\n$,
and thus provides a 'slot' by which we can define a map $S_{(\n)} \rightarrow
E^\sigma_m$ without missing any $M$-degree in $T_\n$. This leads to a somewhat
different philosophy of resolutions than the usual one --- instead of a generating set
of $E^\sigma$ as basic input for our resolutions, we use a polyhedral decomposition.
This at least leads to finite resolutions and reduces in many cases the problem to
understanding the modules $S_{(\n)}$ (see theorem \ref{CMresolution} for such an
application). We want to remark that our notions of admissible posets and polyhedral
decompositions are very close, though not entirely identical, to the sector
partitions in \cite{helmmiller}.

\paragraph{Gluing of posets and globalization.}

A sheaf \msh{E} is equivalent to a collection of $k[\sigma_M]$-modules $E^\sigma$,
where $\sigma$ runs over the fan associated to $X$, which glues in an appropriate
sense over the $U_\sigma$. On the other hand, \msh{E} can be represented by a
collection of sheaves over some admissible posets $\p^\sigma$, which we have to glue
--- in an appropriate sense. The problem of gluing posets might be interesting in a
somewhat broader mathematical context, so that we decided to define it slightly
more general than necessary.
We remark that the naive idea of gluing posets like topological spaces,
which of course can be done, probably does not lead to anything interesting. For
instance, one can easily show that a topological space which is covered by two open
sets, each of which is homeomorphic to a poset, can globally be given the structure of
a poset. By induction, one concludes that every set which has a finite cover by
posets is a poset again.
Our notion of gluing is different from this, and indeed it is a derived
concept which comes very naturally from toric geometry, suitable for us to
construct global resolutions.

Our idea is to realize gluing by passing from posets to {\em preordered} sets.
In contrast to our statements above, a semigroup $\sigma_M$ in general induces only
a preorder on $M$, rather than a partial order. For any two $m, m' \in M$ we have
$m \leq_\sigma m'$ and $m' \leq_\sigma m$ iff $m - m' \in \sigma_M^\bot$, the
maximal subgroup of $\sigma_M$. We can turn $\leq_\sigma$ into a proper partial order
if we pass to the induced order on the quotient $M / \sigma_M^\bot$. For any
$M$-graded module $E$ and any pair $m, m'$ with $m \leq_\sigma m'$ and
$m' \leq_\sigma m$, the multiplication homomorphisms by $\chi(m' - m)$ and
$\chi(m - m')$ necessarily are isomorphisms, and in fact, the categories of
$M$-graded $k[\sigma_M]$-modules and of $M / \sigma_M^\bot$-graded
$k[\sigma_M / \sigma_M^\bot]$-modules are equivalent. Now for simplicity assume that
$\leq_\sigma$ is a partial order and let $\tau < \sigma$ be a proper face, such that
$\leq_\tau$ is a proper preorder. $\tau_M$ is of the form $\sigma_M + \Z_{\geq 0}
\cdot (-m_\tau)$ for some $m_\tau \in \sigma_M$ such that $\tau^\bot \cap
\check{\sigma}$ is a proper face of $\check{\sigma}$ and $\tau^\bot_M =
(\sigma_M \cap \tau_M^\bot) + \Z_{\geq 0} \cdot (-m_\tau)$ is a nontrivial subgroup.

The set $\tau_M^\bot \cap \sigma_M$ is a subsemigroup of $\tau_M^\bot$, giving rise
to a partial order on $\tau_M^\bot$. For any $m \in M$, we can think of the affine
subset $m + \tau_M^\bot$ as a {\em slice} in $M$, and every such slice has its own
partial order. With respect to such a slice, we can consider the directed system
$E^\sigma_{m'}$ with $m' \in m + \tau_M^\bot$, and the directed limit of this system:
\begin{equation*}
E^\tau_m := \underset{\rightarrow}{\lim} E^\sigma_{m'}, \qquad m' \in m + \tau^\bot_M.
\end{equation*}
It turns out that $\bigoplus_{m \in M} E^\tau_m \cong E^\sigma_{\chi(m_\tau)}$, i.e.
the localization of
$E^\sigma$ by the character $\chi(m_\tau)$, which we now can interpret as some kind
of limit figure of $E^\sigma$ along the direction $m_\tau$.

This example is our prototype for defining gluing of partially ordered sets and
sheaves over them. Let $\p$ be an abstract poset with some partial order $\leq$. Then
a {\em localization} of $\leq$ is a preorder $\leq'$, such that $x \leq y$ implies
$x \leq' y$, and $x \leq' y$ implies $x \leq w$ for some element $w$ with $w
\leq' y$ and $y \leq' w$. This is the abstract analogoue of the slicing above, where
the preorder $\leq'$ groups together certain subsets of $\p$. By this definition,
if we pass to the quotient $\p / \sim$, where $x \sim y$ iff $x \leq' y$ and $y
\leq' x$, every representation $F$ of $\p$ induces a representation $\bar{F}$ of
$\p / \sim$, where
\begin{equation*}
\bar{F}_{[x]} = \underset{\rightarrow}{\lim} F_y,
\end{equation*}
the limit is taken over all $y \in x$ {\em with respect to the partial order
$\leq$}. This way, we obtain a quite canonical procedure for gluing sheaves $F_1,
F_2$ over two posets $\p_1, \p_2$; we simply require that the posets have
localizations $\leq'_1, \leq'_2$ such that there exists isomorphisms $l : \p_1 / \sim
\longrightarrow \p_2 / \sim$ and $\phi : l^*\bar{F}_2 \longrightarrow \bar{F}_1$.
We refer to subsections \ref{posetgluing} and \ref{sheafgluing} for the precise
definitions.

\paragraph{Overview of the paper.}

This paper tries to be self-contained in the sense that all required notation related
to toric geometry are introduced. However, we refrain from
giving any account on these subjects; we refer to \cite{perling1} for a more details.

The paper consists of four principal parts. In section \ref{posheaves}, we present
our principal technical framework from the theory of poset representations; in
addition to well-known material, in this section gluing of posets and of sheaves
over posets are introduced.
In section \ref{toricvarieties}, we recall general notions from toric geometry and
the formalism of $\Delta$-families as developed in \cite{perling1}. We present a
partial reformulation of the material in view of the formalism of section
\ref{posheaves}. We show that the Krull-Schmidt theorem holds in the category of
equivariant coherent sheaves over {\em any} toric variety.
Section \ref{resolutions} contains the biggest part of the work;
starting from polynomial rings (subsection \ref{polynomialrings}), and then
generalizing to normal semigroup rings (subsection \ref{semigrouprings}), we construct
resolutions for finitely generated modules over affine toric varieties. In subsections
\ref{homext} and \ref{deltaglobres} we construct global resolutions, both from the
point of view of gluing over posets, and homogeneous coordinate rings.
In section \ref{reflexivesheaves} we analyze the special case of reflexive modules,
and in particular we amplify their close relationship to vector space arrangements.
As an application, in subsection \ref{cmmodules} we show that our resolutions in the
case of reflexive modules behave well in sense of homological algebra. In subsection
\ref{reflexivemodels} we discuss how resolutions of vector space arrangements can
effectively be computed in terms of associated modules over polynomial rings.

\

This work extends results of my thesis \cite{perlingdiss}. Most of this paper
has been written during my stay at the Abdus Salam ICTP, Trieste for whose
hospitality I am deeply grateful.

\section{Preliminaries on Preordered Sets}
\label{posheaves}

In this section let $\mathcal{P}$ be a countable set on which a preorder
$\leq$ is defined. Recall that a preorder is defined by the same axioms as a partial
order, except for the reflexivity axiom, i.e. there may exist elements $x, y \in \p$
such that $x \leq y$ and $y \leq x$, but $x \neq y$. For such pairs we write $x
\lessgtr y$; in the sequel we will frequently put indices on the symbol $\leq$, such
as $\leq', \leq_\sigma$, etc.; then these indices also apply to $\lessgtr$.
If there is no ambiguity in the preorder chosen, we just write $\p$, else we write
$(\p, \leq)$.

\subsection{Representations of preordered sets}

Any preordered set $\mathcal{P}$ in a natural way forms
a category; its objects are given by the set underlying
$\mathcal{P}$ and the morphisms for $x, y \in \ob(\mathcal{P})$ are:
\begin{equation*}
\mor(x, y) =
\begin{cases}
\text{the pair } (x, y) & \text{ if } x \leq y \\
\emptyset & \text{ else.}
\end{cases}
\end{equation*}
Here the pair $(x, x)$ represents the identity morphism for all $x \in \mathcal{P}$.

\begin{definition}
A functor from $\mathcal{P}$ to $\kvect$, the category of vector spaces over the
field $k$, is called a {\em $k$-linear representation} of $\mathcal{P}$.
\end{definition}
As a general notation, if $E$ denotes a $k$-linear representation of a preordered set
$\p$, an element $x \in \p$ is mapped to the vector space denoted $E_x$, and the
relation $x \leq y$ is mapped to a vector space homomorphism $E(x, y)$. 
The $k$-linear representations of $\mathcal{P}$ form an abelian category whose
morphisms are natural transformations.

Representations of $\p$ are equivalent to {\em sheaves} over $\p$. On $\p$ there is
defined a topology which is generated by the basis
\begin{equation*}
U(x) := \{y \geq x\}
\end{equation*}
for all $x \in \mathcal{P}$. The continuous maps between posets then are precisely the
order preserving maps.
A sheaf \msh{E} on $\mathcal{P}$ with respect to this topology with values in $\kvect$
automatically induces a representation of $\mathcal{P}$. On the other hand, for any
representation $E$, following \cite{EGAI} \S 0.3.2,
we obtain a presheaf \msh{E} on $\mathcal{P}$ by setting $\sh{E}\big((U(x)\big)
:= E_x$ for all $x \in \mathcal{P}$ and $\sh{E}(U) := \underset{\leftarrow}{\lim}
\ \sh{E}\big(U(x)\big)$ for some open set $U$, where
the limit runs over all $x \in U$. Note that the stalk $\sh{E}_x$ is isomorphic
to $\sh{E}\big(U(x)\big)$. By observing that for some $U(x)$ every open cover
necessarily contains $U(x)$, and applying the criterion of \S 0.3.2.2 in \cite{EGAI},
it follows that every presheaf automatically is a sheaf.

A distinguished class of representations are the representations $F^x$ associated to
some element $x \in \p$, which are given by:
\begin{equation*}
y \mapsto
\begin{cases}
k & \text{ if } x \leq y \\
0 & \text{ else},
\end{cases}
\end{equation*}
and relations $y \leq z$ mapped to identity if $x \leq y$, and to the zero map
else. In terms of sheaves over $\p$, one can alternatively define $F^x$ as follows.
Denote $j_x$ the canonical inclusion $U(x) \hookrightarrow \p$, and let $\mathbf{k}$
be the constant sheaf with $\mathbf{k}(U(y)) = k$ for all $y \in \p$.
Then $F^x$ corresponds to $j_{x!}j_x^* \mathbf{k}$. We say that a representation
of $\p$ is {\em free} if it is isomorphic to a direct sum of objects of the form
$F^x$. We have:

\begin{proposition}
The representations $F^x$ are projective objects in the category of $k$-linear
representations of $\p$.
\end{proposition}

Using the notion of free objects, we can introduce {\em free resolutions}.

\begin{definition}
Let $\p$ be any preordered set and $E$ a $k$-linear representation.
Then a {\em free resolution} of $E$ is an exact sequence
\begin{equation*}
\dots \longrightarrow F_i \longrightarrow \dots \longrightarrow F_0 \longrightarrow
E \longrightarrow 0
\end{equation*}
where for every $i$ : $F_i \cong \bigoplus_{j} F^{x_{ij}}$ for some $x_{ij} \in \p$.
\end{definition}

Let $x \in \p$, then we consider the subvector space of $E_x$ which is generated
by the image of all $E_y$, $y < x$, by the morphisms $E(x, y)$,
$E_{<x} := \sum_{y < x} E(y, x) E_y$, where we set $E_{<x} := 0$ if the
set $\{y < x\}$ is empty. $\codim_{E_x} E_{< x}$ is
the {\em free dimension} of $E_x$.

\begin{proposition}
\label{repres}
Let $\p$ be a finite preordered set. Then for every $k$-linear representation of
$\p$ there exists a {\em finite} free resolution, that
is, there exists a free resolution as above and some $n \geq 0$ such that $F_i = 0$
for all $i > n$.
\end{proposition}

\begin{proof}
Let $\mathcal{X} \subset \p$ be the set of elements such that $E_x$ has positive free
dimension. For every $x \in \mathcal{X}$ we consider the short exact
sequence of vector spaces
\begin{equation*}
\xymatrix{
0 \ar[r] & E_{<x} \ar[r] & E_x \ar[r] & E_x / E_{<x} \ar[r]
\ar@/^1pc/@{.>}[l]^{\mu_x} & 0,
}
\end{equation*}
where we have chosen some section $\mu_x$. For every such $x$, we can consider
the constant sheaf $E_x / E_{<x}$ on $\p$ and its restriction $E^x := j_{x!}j^*_x (
E_x / E_{<x})$. Using the section $\mu_x$, there exists a natural homomorphism
$\phi_x : E^x \longrightarrow E$ by setting $\phi_x = E(x, y) \circ \mu_x :
E^x_y \longrightarrow E_y$
for every pair $x \leq y$ and the zero map for all $x \nleq y$.
The sheaf $E^x$ is isomorphic to $(F^x)^{f_x}$, where $f_x$ is the free dimension of
$E_x$. Thus we define $F_0 = \bigoplus_{x \in \mathcal{X}} E^x \cong
\bigoplus_{x \in \mathcal{X}} (F^x)^{f_x}$ and a homomorphism
$\phi_0: F_0 \longrightarrow E$ by setting $\phi_0 := \sum_{x \in \mathcal{X}}
\phi_x$. By construction, $\phi_0$ is a surjective map, and we obtain thus
a short exact sequence of representations of $\p$:
\begin{equation*}
0 \longrightarrow K_0 \longrightarrow F_0 \overset{\phi_0}{\longrightarrow} E
\longrightarrow 0.
\end{equation*}
Now we can repeat this construction with $K_0$, and by iterating we obtain a free
resolution of $E$ which is concatenated of short exact sequences $0 \longrightarrow
K_{i + 1} \longrightarrow F_{i + 1} .\overset{\phi_{i + 1}}{\longrightarrow} K_i
\longrightarrow 0$.

Now observe that for $K_{i + 1, x} = 0$ whenever the free dimension of $K_{i, x}$ is
equal to $\dim K_{i, x}$,
and $K_{i, x} = 0$ implies that $K_{i + 1, x} = 0$. The set of such 
$K_{i, x}$ whose free dimension is equal to $\dim K_{i, x}$ is always nonempty as long
as $K_i$ is nontrivial, because the set contains at least the minimal elements
$x \in \p$ which have nontrivial $K_{i, x}$. So, as $\p$ is finite, it follows that
there exists some $r > 0$ for which $K_{i, x} = 0$ for all $i > r$.
\end{proof}

\begin{definition}
Let $0 \longrightarrow F_r \overset{\phi_r}{\longrightarrow} \dots,
\overset{\phi_1}{\longrightarrow} F_0 \overset{\phi_0}{\longrightarrow} E
\longrightarrow 0$ be a free resolution of a $k$-linear representation $E$, then we
call the kernel of $\phi_i$ the {\em $i$th syzygy representation} of $E$.
\end{definition}

\subsection{Direct and inverse limits}

Now we recall some basic facts about direct and inverse limits in the
category of vector spaces. This is only intended as a reminder to the reader,
as we will be using limits extensively during the rest of this paper.

As we have seen in the previous subsection, every preordered set $\p$ in a natural
way is a directed family. Thus, a representation $E$ of $\p$ becomes a directed
family of vector spaces. Recall, that the {\em inverse limit} of $E$ is a vector
space
\begin{equation*}
\underset{\leftarrow}{\lim} E =: \mathbf{E}^i
\end{equation*}
which has the following universal properties:
\begin{enumerate}[(i)]
\item for every element $x \in \p$ there exists a unique homomorphism $\phi_x :
\mathbf{E}^i \longrightarrow E_x$ such that $E(x, y) \circ \phi_x = \phi_y$ for
every $x \leq y$;
\item for every vector space $\mathbf{F}$ with homomorphisms $\psi_x: \mathbf{F}
\longrightarrow E_x$, where $\psi_y = E(x, y) \circ \psi_x$ for every $x \leq y$,
there exists a unique homomorphism $\delta: \mathbf{F} \longrightarrow \mathbf{E}^i$
with $\psi_x = \phi_x \circ \delta$ for all $x \in \p$.
\end{enumerate}

\begin{definition}
We denote the vector space homomorphism $\delta: \mathbf{F} \longrightarrow
\mathbf{E}^i$ {\em diagonal homomorphism} from $\mathbf{F}$ to $\mathbf{E}^i$.
\end{definition}

Explicitly, such a limit can be constructed as the subvector space of the direct
product $\prod_{x \in \p} E_x$ consisting of sequences $(e_x \mid x \in \p)$ such
that $E(x, y)(e_x) = e_y$ for every pair $x \leq y$. If $\p$ has a unique minimal
element $x_{\min}$, then $\phi_{x_{\min}} : \mathbf{E}^i \rightarrow E_{x_{\min}}$
becomes an isomorphism. This construction is a straightforward
generalization of the
pullback in the category of vector spaces; the pullback is the special case where
the poset consists of three elements $x, y, z$ with $x < z$ and $y < z$.

Dually, there exists the {\em direct limit}
\begin{equation*}
\underset{\rightarrow}{\lim} E =: \mathbf{E}^d
\end{equation*}
which generalizes pushout. It can explicitly be constructed as the quotient of
the vector space $\prod_{x \in \p} E_x$ by the subvector space generated by vectors
$E(x, y)(e_x) - e_x$. For every $x \in \p$ there exists a homomorphism $\phi^x :
E_x \longrightarrow \mathbf{E}^d$ such that universal properties analogously to the
inverse limit are fulfilled. Note that in case that there exists a unique maximal
element $x_{\max}$, the homomorphism $\phi^{x_{\max}}$ is an isomorphism.

Both limits behave covariantly; consider two preordered sets $\p$, $\mathcal{Q}$,
and any two representations $E$, $F$ of $\p$ and $\mathcal{Q}$, respectively, and
a order preserving map $f : \p \longrightarrow \mathcal{Q}$. Then any natural
transformation $r : E \longrightarrow f^* F$ induces a homomorphism of limits:
\begin{equation*}
\underset{\leftarrow}{\lim}\ r : \underset{\leftarrow}{\lim}\ E  \longrightarrow
\underset{\leftarrow}{\lim}\ f^*F \text{\quad resp. \quad }
\underset{\rightarrow}{\lim} \ r : \underset{\rightarrow}{\lim}\ E  \longrightarrow
\underset{\rightarrow}{\lim}\ f^*F.
\end{equation*}
In particular, if $\p$ and $\mathcal{Q}$ have unique minimal elements $x_{\min}$ and
$y_{\min}$, respectively, and $f(x_{\min}) = y_{\min}$, we obtain
\begin{equation*}
\underset{\leftarrow}{\lim}\ r : \underset{\leftarrow}{\lim}\ E  \longrightarrow
\underset{\leftarrow}{\lim}\ F,
\end{equation*}
and analogously for the direct limit with respect to maximal elements.

\subsection{Gluing of preordered sets}
\label{posetgluing}

\begin{definition}
Let $(\p, \leq)$ be a preordered set. Then we denote $\p_\lessgtr$ the quotient
of $\p$ by the equivalence relation which is given by $x \sim y$ iff $x \lessgtr y$.
\end{definition}

Clearly, $\leq$ induces a partial order on the set $\p_\lessgtr$.

\begin{definition}
A {\em localization} of $\p$ is a preorder $\leq'$ on $\mathcal{P}$ such that the
following conditions are fulfilled:
\begin{enumerate}[(i)]
\item for all $x, y \in \p$, $x \leq y$ implies $x \leq' y$,
\item for all $x \leq' y$ there exists some $w \lessgtr' y$ such that
$x \leq w$.
\end{enumerate}
\end{definition}

Let $\leq'$ be a localization of $(\p, \leq)$, then $\leq$ induces a relation on
$\p_{\lessgtr'}$ by setting $[x] \leq [y]$ iff there exist $u \in [x]$, $v \in [y]$
with $u \leq v$.

\begin{proposition}
Let $\leq'$ a localization of $(\p, \leq)$, then the relation on $\p_{\lessgtr'}$
induced by $\leq$ coincides with the partial order induced by $\leq'$.
\end{proposition}

\begin{proof}
We check the poset axioms for $\leq$: {\em 1)} $[x] \leq [x]$ follows because
$u \leq' u$ implies that there exists $v \lessgtr' u$ such that $u \leq v$. {\em 2)}
Let $[x] \leq [y]$ and $[y] \leq [x]$; then there exist $u, p \in [x]$, $v, q \in
[y]$ such that $u \leq v$ and $p \leq q$; then $u \leq'v \leq' p \leq' u$, and thus
$v \lessgtr' u$, hence $[x] = [y]$. {\em 3)} Let $[x] \leq [y]$ and $[y] \leq [z]$;
then there exist $u \in [x]$, $v, p \in [y]$ and $q \in [z]$ such that $u \leq v$ and
$p \leq q$; thus $u \leq' v \leq' p \leq' q$ and there exists $w \lessgtr'q$ such
that $u \leq w$, and thus $[x] \leq [z]$.

Now the equivalence of the partial orders $\leq$ and $\leq'$ on $\p_{\lessgtr'}$ is
trival.
\end{proof}

Let $(\p_\alpha, \leq_\alpha)$ be a finite family of preordered sets, together
with a family of representations $F_\alpha$, where $\alpha$ runs over some index set
$A$. Our aim is to {\em glue} these representations when certain conditions on the
$\p_\alpha$ are fulfilled. For this, we need the following notion:

\begin{definition}
Let $f : \p \longrightarrow \mathcal{Q}$ be an order preserving
map between preordered sets. Then $f$ is a {\em contraction} if
\begin{enumerate}[(i)]
\item for every $x \in \p$ there exists some $y \in \mathcal{Q}$ such that
$f\big(U(x)\big) = U(y)$,
\item for every $y \in \mathcal{Q}$: $f^{-1}\big(U(y)\big) = U(x)$ for some $x \in
\p$.
\end{enumerate}
\end{definition}
These conditions imply that $f$ is surjective and that for every $x \in
\p$ with $f\big(U(x)\big) = U(y)$ there exists $z \in \p$ with $U(x) \subset U(z) =
f^{-1}\big(U(y)\big)$. By this we can define a map $h : \mathcal{Q} \longrightarrow
\p$ by mapping $y \mapsto z$. This map is an order preserving injection of
$\mathcal{Q}$ into $\p$.

\begin{definition}
Let $f : \p \longrightarrow \mathcal{Q}$ be a contraction. Then the unique map $h :
\mathcal{Q} \longrightarrow \p$ mapping $y \in \mathcal{Q}$ to $z \in \p$ such that
$f^{-1}(U(y)) = U(z)$ is called {\em hooking} of $\mathcal{Q}$ into $\p$.
\end{definition}

Using this definition, one can think of our gluing of posets as a process of
{\em hooking} different posets along common contractions.

Let $\p, \mathcal{Q}$ be two finite preordered sets, $f : \p \longrightarrow
\mathcal{Q}$ a contraction, and $E$ a representation of $\mathcal{Q}$. Then the
pullback $f^* F^x$ for any free representation for some $x \in \mathcal{Q}$ then
is isomorphic to the free representation $F^{h(x)}$ of $\p$. For any free
resolution $0 \rightarrow F_r \rightarrow \cdots \rightarrow F_0 \rightarrow E
\rightarrow 0$, one can consider the pullback sequence $0 \rightarrow f^*F_r
\rightarrow \cdots \rightarrow f^*F_0 \rightarrow f^*E \rightarrow 0$. We observe:
\begin{lemma}
\label{contractionliftres}
The sequence $0 \rightarrow f^*F_r \rightarrow \cdots \rightarrow f^*F_0 \rightarrow
f^*E \rightarrow 0$ is isomorphic to the free resolution of $f^*E$ in the sense of
proposition \ref{repres}.
\end{lemma}

\begin{proof}
It suffices to check the first step of the resolution $0 \rightarrow K_0 \rightarrow
f^* F_0 \rightarrow f^* E \rightarrow 0$ and to show that $K_0$ and $f^* F_0$
coincide with the representations obtained by the procedure of proposition
\ref{repres}. But this follows directly from the fact that for every $y \in
\mathcal{Q}$, the homomorphisms $(f^*E)_{h(x)} \rightarrow (f^*E)_y$ are isomorphisms
for all $h(x) \leq y \in f^{-1}(x)$.
\end{proof}

\begin{definition}
Let $(A, \preceq)$ be a finite poset and $\p_\alpha$, $\alpha \in A$ be a
family of preordered sets. We say that the posets $\p_\alpha$ {\em glue over $A$} if
\begin{enumerate}[(i)]
\item for every $\beta < \alpha \in A$ there exists a localization
$\leq_\alpha^\beta$ of $\leq_\alpha$ and a contraction 
$l_{\alpha\beta}: (\p_\alpha)_{\lessgtr_\alpha^\beta} \rightarrow
(\p_\beta)_{\lessgtr_\beta}$;
\item for every triple $\gamma
\preceq \beta \preceq \alpha \in A$, the composition of maps $\p_\alpha
\rightarrow (\p_\alpha)_{\lessgtr^\beta_\alpha}
\overset{l_{\alpha\beta}}{\rightarrow} (\p_\beta)_{\lessgtr_\beta}
\rightarrow (\p_\beta)_{\lessgtr_\beta^\gamma}\overset{l_{\beta\gamma}}{\rightarrow}
(\p_\gamma)_{\lessgtr_\gamma}$ coincides with $\p_\alpha \rightarrow
(\p_\alpha)_{\lessgtr^\gamma_\alpha}
\overset{l_{\alpha\gamma}}{\rightarrow} (\p_\gamma)_{\lessgtr_\gamma}$.
\end{enumerate}
\end{definition}

Our principal example, where the maps $l_{\alpha\beta}$ actually are isomorphisms,
will be the preorderings associated to a fan $\Delta$ in section
\ref{toricvarieties}.

\subsection{Gluing of sheaves over preordered sets}
\label{sheafgluing}

Let $(\p, \leq)$ be some preordered set; if $E$ is some representation of $\p$, then
for any pair $x \lessgtr y$, the map $E(x,y) : E_x \longrightarrow E_y$ is an
isomorphism whose inverse is $E(y,x)$. Thus $E$ descends to a representation of
$\p_\lessgtr$ by setting $E_{[x]} := \underset{\rightarrow}{\lim} E_y$, where
the direct limit is taken over all elements $y \leq x$. 
For any $y \leq x$ there is the canonical inclusion of directed systems $\{E_z \mid
z \leq y\} \hookrightarrow \{E_z \mid z \leq x\}$, which induces a functorial
homomorphism
$E_{[y]} \longrightarrow E_{[x]}$. On the other hand, every representation $F$ of
$\p_\lessgtr$ lifts to a representation of $\p$ by setting $E_x := E_{[x]}$ and
$E(x,y) := E([x],[y])$. By descend and lift, we have:

\begin{lemma}
\label{popreequivrep}
Let $(\p, \leq)$ be a preordered set. The category of representations of $\p$ is
equivalent to the category of representations of $\p_\lessgtr$.
\end{lemma}

\comment{
Let $(\p, \leq_1)$ be a preorderet set and let $\leq_2$ be another preorder on
$\p$ such that $x \leq_1 y$ implies $x \leq_2 y$ for all $x, y \in \p$. Denote
$\sim_i$ the equivalence relation with respect to $\leq_i$, $i = 1, 2$; the
preorder $\leq_2$ induces a preorder on $\p / \sim_1$, and forming equivalence
classes of $\p / \sim_1$ with respect to the induced preorder coincides with
$\p / \sim_2$.
}

Let $\leq'$ be a localization of $\leq$. For any $x \in \p$, denote
$\p_x = \{z \in \p \mid z \leq' x\}$. We construct a representation on
$\p_{\lessgtr'}$ by mapping $[x]' \in (\p)_{\lessgtr'}$ to the vector space $E_{[x]'}
:= \underset{\rightarrow}{\lim} E_z$, the direct limit taken over
$\p_x$ {\em with respect to the partial order $\leq$}.
The inclusion $\p_x \hookrightarrow \p_y$ induces an inclusion of directed sets
with respect to $\leq$, and thus we obtain a morphism $E_{[x]'} \longrightarrow
E_{[y]'}$. By lemma \ref{popreequivrep}, this representation lifts to
a representation of $(\p, \leq')$.

\begin{definition}
Let $(\p, \leq)$ be a preordered set, $\leq'$ a localization of $\leq$, and $E$ a
representation of $(\p, \leq)$. Consider the poset $\p_{\lessgtr'}$. Then we call
the induced sheaf on $(\p, \leq')$ a {\em localization} of $F$.
\end{definition}

Now we assume that we are given some partially ordered set $(A, \preceq)$, a
collection of preordered sets $\p_\alpha$ which glues over $A$, and a collection of
sheaves $E^\alpha$ over $\p_\alpha$ for every $\alpha \in A$. We want {\em glue} this
collection of sheaves to give some kind of global object over the glued preordered
sets.

\begin{definition}
We say that the collection $E^\alpha$ {\em glues} over the collection $\p_\alpha$, if
\begin{enumerate}[(i)]
\item for every $\beta \preceq \alpha$, and morphism of posets $l_{\alpha\beta}:
(\p_\alpha)_{\lessgtr^\beta_\alpha} \longrightarrow (\p_\beta)_{\lessgtr_\beta}$
there is an isomorphism of sheaves $\phi^{\alpha\beta} : l_{\alpha\beta}^*E^\beta
\overset{\cong}{\longrightarrow} E^\alpha$
\item for every triple $\gamma \preceq \beta \preceq \alpha$: $\phi^{\alpha\gamma} =
\phi^{\alpha\beta} \circ l^*_{\alpha\beta}\phi^{\beta\gamma}$.
\end{enumerate}
We call a such a collection a {\em sheaf} over $\p^\alpha$.
\end{definition}

Let $E^\alpha$, $F^\alpha$ be sheaves over $\p^\alpha$, where we denote the gluing
homomorphisms $\phi^{\alpha\beta}$ and $\psi^{\alpha\beta}$, respectively.
A homomorphism from $E^\alpha$ to $F^\alpha$ is given by a collection of
homomorphisms $f_\alpha : E^\alpha \longrightarrow F^\alpha$ such that $f_\alpha
\circ \phi^{\alpha\beta} = \psi^{\alpha\beta} \circ l^*_{\alpha\beta}f_\alpha$
for every pair $\beta \preceq \alpha$. One checks straightforwardly that this is
compatible with the cocycle conditions on $\phi^{\alpha\beta}$ and
$\psi^{\alpha\beta}$, and moreover that the corresponding families of kernels and
cokernels of $f^\alpha$ glues over $A$:

\begin{proposition}
The category of sheaves over $\p^\alpha$ is abelian.
\end{proposition}

\comment{
Now we define an anologon for invertible sheaves:

\begin{definition}
Let $\p^\alpha$ be a collection of posets which glues over $A$ and consider some
collection $\{x_\alpha \in \p_\alpha\}$ which has the property that for every
$\beta < \alpha$ $\bar{x}_\alpha = h_{\alpha\beta}(x_\beta)$. Then a {\em locally
free} representation of the system $x_\alpha$ is given by a collection $F^{x_\alpha}$
which glues over the collection $\p_\alpha$.
\end{definition}
}

\paragraph{Compression of sheaves over preordered sets.}

Let $A$ be a finite poset and denote $\mathbf{P}^f_A$ the category of collections of
{\em finite} preordered sets $\{\p^\alpha \mid \alpha \in A\}$ which glue over $A$.
Let $\{\mathcal{Q}^\alpha\}$ be any collection of not necessarily finite preordered
sets which glues over $A$. Denote $\mathbf{C}$ any subcategory of the category of
sheaves which glue over the collection $\mathcal{Q}^\alpha$. A {\em compression} of
$\mathbf{C}$ is any object $\{\p^\alpha\}$ of $\mathbf{P}^f$ together with a pair of
functors
\begin{align*}
\zip & : \mathbf{C} \longrightarrow \mathbf{Sheaves}(\p^\alpha) \\
\unzip & : \mathbf{Sheaves}(\p^\alpha) \longrightarrow \mathbf{C}.
\end{align*}
which induce an equivalence of categories between $\mathbf{C}$ and
$\mathbf{Sheaves}(\p^\alpha)$.

\section{Toric Varieties and $\Delta$-Families}
\label{toricvarieties}

In this section we briefly recall basic facts for toric varieties and our results
from \cite{perling1} on equivariant sheaves over toric varieties.
For general information about toric varieties we refer
to \cite{Oda} and \cite{Fulton}. In this work $X$ will always denote an
$r$-dimensional
toric variety over a fixed algebraically closed field $k$, and $T$ the open dense
torus contained in $X$. Moreover, we use the following notation:

\begin{itemize}
\setlength{\itemsep}{-4pt}
\item $M \cong \mathbb{Z}^n$ is the character group of $T$, and $N$ the
$\mathbb{Z}$-module dual to $M$,\\
$M_\R := M \otimes_Z \R$, $N_\R := N \otimes_Z \R$;
\item elements of $M$ are denoted $m, m'$ etc. if written additively and $\chi(m),
\chi(m')$ etc. if written multiplicatively, i.e. $\chi(m + m') = \chi(m)\chi(m')$;
\item $\Delta$ denotes the fan associated to $X$, and cones in $\Delta$ are denoted
by small Greek letters $\rho$, $\sigma$, $\tau$, etc.;
the natural order among cones is denoted by $\tau < \sigma$, \\
$\Delta(i) := \{\sigma \in \Delta \mid \dim \sigma = i\}$ the set of
all cones of fixed dimension $i$,\\
$\sigma(i) := \{\tau \in \Delta(i) \mid \tau < \sigma\}$;
\item $\check{\sigma} := \{m \in M_\mathbb{R} \mid \langle m, n \rangle \geq 0
\text{ for all $n \in \sigma$}\}$ is the cone {\it dual} to $\sigma$,\\
$\sigma^\bot = \{m \in M_\mathbb{R} \mid \langle m, n \rangle = 0 \text{ for all } n \in \sigma \}$, \\
$\sigma_M := \check{\sigma} \cap M$ is the subsemigroup of $M$ associated to
$\sigma$. \\
$\sigma_M^\bot := \sigma^\bot \cap M$ is the maximal subgroup of $\sigma_M$;
\item the affine toric variety associate to a cone $\sigma$ is denoted $U_\sigma$,\\
$U_\sigma \cong \spec{k[\sigma_M]}$, where $k[\sigma_M]$ is the semigroup ring over
$\sigma_M$;
\item elements of $\Delta(1)$ are called {\it rays}, and the torus invariant Weil
divisor associated to some ray $\rho \in \rays$ is denoted $D_\rho$.
\end{itemize}

\subsection{Equivariant sheaves and $\Delta$-families}
\label{deltafamilies}

Consider any rational polyhedral convex cone $\sigma$, then the subsemigroup
$\sigma_M$ induces a {\em directed preorder} $\leq_\sigma$ on $M$ by setting
$m \leq_\sigma m'$ iff $m' - m \in \sigma_M$. The following properties of
$\leq_\sigma$ are easy to see:
\begin{enumerate}[(i)]
\setlength{\itemsep}{-5pt}
\item $m \leq_\sigma m'$ and $m' \leq_\sigma m$ iff $m - m' \in \sigma_M^\bot$.
\item If $\tau \leq \sigma$, then $m \leq_\sigma m'$ implies $m \leq_\tau m'$.
\item If $\sigma$ is of maximal dimension in $N_\R$, then $\leq_\sigma$ is a partial
order.
\end{enumerate}

Let $\sh{E}$ be an equivariant sheaf over $X$ and denote $E^\sigma :=
\Gamma(U_\sigma, \sh{E})$ for every affine open $T$-invariant subvariety $U_\sigma$
of $X$. The dual action of $T$ on $E^\sigma$ induces an isotypical decomposition
\begin{equation*}
E^\sigma = \bigoplus_{m \in M}E^\sigma_m
\end{equation*}
For any two $m \leq_\sigma m'$, there exists a distinguished $k$-linear map
spaces $\chi^\sigma_{m, m'} : E_m \longrightarrow E_{m'}$ which is given by
multiplication by the monomial
$\chi(m' - m) \in k[\sigma_M]$. These distinguished maps completely specifiy the
module structure of $E^\sigma$ over $k[\sigma_M]$. Observing that $\chi(m'' - m')
\chi(m' - m) = \chi(m'' - m)$ and $\chi(m - m) = 1$, we even obtain a
{\em functorial} description of $E^\sigma$. By mapping $m \mapsto E^\sigma_m$ and
$(m, m') \mapsto \chi^\sigma_{m, m'}$ for $m \leq_\sigma m'$, every $M$-graded
$k[\sigma_M]$-module
$E^\sigma$ defines a functor from the preordered set $(M, \leq_\sigma)$ to the
category $\kvect$ of $k$-vector spaces.

\begin{proposition}[\cite{perling1}, Proposition 5.5]
Let $U_\sigma = \spec{k[\sigma_M]}$ be an affine toric variety. Then the following
categories are equivalent:
\begin{enumerate}[(i)]
\setlength{\itemsep}{-4pt}
\item equivariant quasicoherent sheaves over $U_\sigma$,
\item $M$-graded $k[\sigma_M]$-modules,
\item $k$-linear representations of the preordered set $(M, \leq_\sigma)$.
\end{enumerate}
\end{proposition}

\begin{definition}
We call a representation of $(M, \leq_\sigma)$ a {\em $\sigma$-family}.
\end{definition}

In the sequel, we will use the notation $E^\sigma$ exchangeably for the
$k[\sigma_M]$-module and for the $\sigma$-family.

Now for any pair $\tau < \sigma$, there exists some $m_\tau \in \sigma_M^\bot$ such
that $\tau_M = \sigma_M + \Z_{\geq 0} (-m_\tau)$ and $\tau_M^\bot = (\tau_M^\bot \cap
\sigma_M) + \Z_{\geq 0} (-m_\tau)$. In terms of preordered sets, this translates the
way that we can consider $(M, \leq_\tau)$ as a localization of $(M, \leq_\sigma)$ in
the sense of subsection \ref{posetgluing}. Moreover, the localization of $(M,
\leq_\sigma)$ by $\leq_\tau$ coincides with $(M, \leq_\tau)$, and thus the
contractions
$l_{\sigma\tau} : M_{\lessgtr_\sigma^\tau} \longrightarrow M_{\lessgtr_\tau}$ are
isomorphisms. We have:

\begin{proposition}
The family of preordered sets $(M, \leq_\sigma)$, $\sigma \in \Delta$, glues over
$\Delta$.
\end{proposition}

The restriction of $\sh{E}\vert_{U_\sigma}$ to $U_\tau$ corresponds to the
localization $E^\sigma_{\chi(m_\tau)}$. To
understand this in terms of $\sigma$-families, we first observe that the canonical
map $E^\sigma \longrightarrow E^\sigma_{\chi(m_\tau)}$ at the same time is a
homomorphism of directed systems.

\begin{proposition}
For every $m \in M$ there exists a natural isomorphism $E^\tau_m \cong
\underset{\rightarrow}{\lim} E^\sigma_{m'}$, where the limit is taken over the
directed
system of all $E^\sigma_{m'}$ with $m' \leq_\tau m$ {\em with respect to the preorder
$\leq_\sigma$}.
\end{proposition}

\begin{proof}
By definition of localization, the vector space $E^\tau_m$ is the set of equivalence
classes $\{[\frac{e}{\chi(m')}] \mid \deg_M e = m + m'\}$, where
$\frac{e_1}{\chi(m_1)}
\sim \frac{e_2}{\chi(m_2)}$ if and only if $\chi(m_1) \cdot e_2 = \chi(m_2) \cdot
e_1$ in $E^\sigma$, where without loss of generality, $m_1$ and $m_2$ can be chosen
from $\sigma_M$. In other notation, this reads $\chi^\sigma_{m + m_1, m + m_1 +
m_2} e_2 = \chi^\sigma_{m + m_2, m + m_1 + m_2} e_1$. So, in a natural way, we can
identify $E^\tau_m$ with the direct limit $\underset{\rightarrow}{\lim}
E^\sigma_{m'}$.
\end{proof}

\comment{ Now for any two $m_1 \lessgtr_\tau
m_2$, the map $\chi^\tau_{m_1, m_2}$ is an isomorphism, and moreover there exist
$m'_1, m'_2 \in \sigma^\bot_M$ such that $\chi(m'_1) \cdot E_{m_1} = \chi(m'_2)
\cdot E_{m_2}$. This in particular implies that there is an isomorphism
\begin{equation*}
E^\tau_m \cong \underset{\rightarrow}{\lim} E^\sigma_{m'}
\end{equation*}
for every $m \in M$, where the limit is taken over the directed system of all
$E^\sigma_{m'}$ with $m' \leq_\tau m$ {\em with respect to the preorder
$\leq_\sigma$}. The canonical map
\begin{equation*}
E^\sigma \longrightarrow E^\sigma_{\chi(m_\tau)}
\end{equation*}
in terms of the $\sigma$- and $\tau$-families translates into a localization of
the representation $E^\sigma$ of $(M, \leq_\sigma)$ to $\leq_\tau$.
}

By this proposition, we see that the localization of $E^\sigma$ by $\chi(m_\tau)$
translates into the localization of $E^\sigma$, considered as {\em sheaf over
$(M, \leq_\sigma)$}, to $(M, \leq_\tau)$. We get:

\begin{definition}[see also \cite{perling1}, Definition 5.8]
A $\Delta$-family is a collection $\{E^\sigma \mid \sigma \in \Delta\}$ of
$\sigma$-families which glues over $\Delta$.
\end{definition}

\begin{theorem}[\cite{perling1}, Theorem 5.9]
The category of equivariant sheaves over $X$ is equivalent to the category of
$\Delta$-families.
\end{theorem}

\subsection{The Krull-Schmidt property}
\label{krullschmidt}

Let $\mathfrak{C}$ be any category in which direct
sums exist. We say that the Krull-Schmidt theorem holds in $\mathfrak{C}$ if
for every object $A$ in $\mathcal{C}$ and for every two decompositions into
indecomposable objects
\begin{equation*}
A \cong X_1 \oplus X_2 \oplus \dots \oplus X_n \cong Y_1 \oplus Y_2 \oplus \dots
\oplus Y_m
\end{equation*}
we have $m = n$, and there exists a permutation $\pi$ of $\{1, \dots, n\}$ such that
$X_i \cong Y_{\pi(i)}$ for every $i$.
It is well known that the Krull-Schmidt theorem holds in the category of coherent
sheaves over a complete variety. For the category of equivariant coherent sheaves
over a toric variety, we can drop the completeness condition:

\begin{theorem}
\label{krullschmidttheorem}
Let $X$ be any toric variety, then the Krull-Schmidt theorem holds for the category
of equivariant coherent sheaves over $X$.
\end{theorem}

\begin{proof}
According to a classical result of Atiyah (\cite{Atiyah1}),
it suffices to show that for every
two equivariant sheaves \msh{E} and $\sh{F}$, the vector space
$\Hom(\sh{E}, \sh{F})^T$ of $T$-equivariant sheaf homomorphisms is
finite-dimensional. As we are dealing only with finite fans, it is enough to consider
the case where $X = U_\sigma$ is an affine toric variety such that $\sh{E}$ and
$\sh{F}$ correspond to finitely generated $k[\sigma_M]$-modules $E^\sigma$ and
$F^\sigma$. In this case the statement follows because every generator of $E^\sigma$
of degree $m$ must be mapped to some element $f \in F^\sigma_m$ and every vector
space $F^\sigma_m$ is finite dimensional (\cite{perling1}, Proposition 5.11).
\end{proof}

\subsection{The quotient representation of a toric variety}

Every toric variety can be represented as a good quotient of a quasi-affine toric
variety (see \cite{Cox}). This representation starts with the exact sequence
\begin{equation*}
0 \longrightarrow M_0 \longrightarrow M \longrightarrow \weildivisors
\longrightarrow A \longrightarrow 0
\end{equation*}
where the map from $M$ to \mweildivisors\ is given by $m \mapsto (\langle m, n(\rho)
\rangle \mid \rho \in \rays)$. In the sequel we will assume that the fan $\Delta$ is
not contained in a proper subvector space of $N_\R$. In this case $M_0$ is the zero
module.

We consider the polynomial ring $S = k[x_\rho \mid \rho \in \rays]$; this
ring is endowed with a natural \mweildivisors-grading by setting $\deg x^\n = \n$ for
every monomial $x_\rho$. Via the surjection of \mweildivisors\ onto $A$, the ring $S$
automatically acquires an $A$-grading,
\begin{equation*}
S \cong \bigoplus_{\alpha \in A} S_\alpha.
\end{equation*}

We define the {\em irrelevant} ideal $B = \langle x^{\hat{\sigma}} \mid \sigma \in
\Delta \rangle$, where $x^{\hat{\sigma}} = \prod_{\rho \in \rays \setminus \sigma(1)}
x_\rho$ for every $\sigma \in \Delta$.
The variety $\mathbf{V}(B)$ defined by $B$ is a finite union of linear subspaces of
$\spec{S} \cong k^\rays$, which has codimension at least two.
The complement of $\mathbf{V}(B)$, which we denote $\hat{X}$, is a quasi-affine toric
variety, on which the torus $\hat{T} \cong (k^*)^\rays$ acts. Denote $e_\rho$ the
standard basis vectors of $\R^\rays$, then the fan of $\hat{X}$ is generated by
the cones $\hat{\sigma} = \sum_{\rho \in \sigma(1)} \R_{\geq 0} e_\rho$, for every
$\sigma \in \Delta$. The affine open subsets
$U_{\hat{\sigma}}$ form a cover of $\hat{X}$, and we will call $\hat{\Delta}
= \{\hat{\sigma} \mid \sigma \in \Delta\}$
the fan of $\hat{X}$, although in general $\hat{\Delta}$ is not a proper fan, unless
$X$ is a simplicial toric variety. 
There is a canonical morphism $\pi: \hat{X} \longrightarrow X$ which is described by
the map of fans induced by the linear map given by $e_\rho \mapsto n(\rho)$. By this
morphism, $X$ becomes a {\em good quotient} of $\hat{X}$ by the diagonalizable
group $G = \Hom(A, k^*)$.
The coordinates $x_\rho$ then serve as global coordinates for 
$X$, and $S$ is denoted the {\em homogeneous coordinate ring}
of $X$.

\comment{
As every $U_\sigma$ is a toric variety, every such $U_\sigma$ has its own homogeneous
coordinate ring $S_\sigma$. Below it will be useful for us to observe that we can
obtain $S$ as a limit ring of the $S_\sigma$.

 This quotient
has locally the description of $U_\sigma = U_{\hat{\sigma}} // G$, and $k[\sigma_M]
= S_{x^{\hat{\sigma}}}^G = \big(S_{x^{\hat{\sigma}}}\big)_0$ with respect to the
$A$-grading.
Note that for any $\sigma \in \Delta$, the localization $S_{x^{\hat{\sigma}}}$
automatically becomes a homogeneous coordinate ring for $U_\sigma$, and the
restriction $\pi \vert_{U_{x^{\hat{\sigma}}}} \longrightarrow U_\sigma$ becomes a
quotient representation for $U_\sigma$.
}

\paragraph{$A$-graded $S$-modules.}
Any $A$-graded $S$-module $F$ defines a $G$-equivariant sheaf over
$k^\rays$ and thus over $\hat{X}$, and it has been shown (see \cite{Mustata1}) that
every quasicoherent sheaf over $X$ can be represented as a descend of an $A$-graded
$S$-module $F$ of the
form $\big(\pi_* (\tilde{F}\vert_{\hat{X}})\big)^G$, where $\tilde{\ }$ denotes the
usual sheafification functor over the affine space $k^\rays$. We abbreviate the
descend of a module $F$ by $\breve{F}$. In the other direction, every coherent sheaf
\msh{F} over $X$ gives rise to an $A$-graded $S$-module $\Gamma(\hat{X},
\pi^*\sh{F})$. There is always an isomorphism $\Gamma(\hat{X}, \pi^*\sh{F})\breve{\ }
\cong \sh{F}$, but in general there is no isomorphism between any $A$-graded module
$F$ and $\Gamma(\hat{X}, \pi^* \breve{F})$.

\paragraph{Fine-graded $S$-modules.}
For the study of equivariant sheaves, we have to
consider {\em fine graded} modules, i.e. \mweildivisors-graded $S$-modules. Such a
module $F$ is equivalent to $\hat{T}$-equivariant sheaf over $k^\rays$, and its
descend $\breve{F}$ then in a natural way is a $T$-equivariant sheaf over $X$.
On the other hand, the pullback $\pi^*\sh{E}$ of some $T$-equivariant sheaf over $X$
has a natural $\hat{T}$-equivariant structure, and thus $\hat{E} := \Gamma(\hat{X},
\pi^* \sh{E})$ is fine graded. The most important examples for us are
the modules which are defined as the descend of free $S$-modules of rank one. These
are the modules of the form $S(\n)$, the degree shift of $S$ by some element
$\n \in \weildivisors$, where $S(\n)_{\n'} = S_{\n + \n'}$. The descend
$\breve{S}(\n)$ is isomorphic to $\sh{O}_X(D_\n)$, the reflexive sheaf of rank one
which is associated to the Weil-divisor $D_\n := \sum_{\rho \in \rays} -n_\rho
D_\rho$. As a general notation, we write $S_{(\n)}$ instead of $S(\n)_0$; note that
this shift is in the \mweildivisors-grading, not in the $A$-grading and therefore
fixes a unique equivariant structure on $\breve{S}(\n)$.

\paragraph{Global and local quotient representations.}
For any $\sigma \in \Delta$ there is an exact sequence
\begin{equation*}
0 \longrightarrow \sigma_M^\bot \longrightarrow M \longrightarrow \Z^{\sigma(1)}
\longrightarrow A^\sigma \longrightarrow 0,
\end{equation*}
by which we have a splitting $M \cong \sigma_M^\bot \oplus M / \sigma_M^\bot$, where
we identify $M / \sigma_M^\bot \cong M_{\lessgtr_\sigma}$ with the image of $M$ in
$\Z^{\sigma(1)}$. This
induces a splitting $U_\sigma \cong T_\sigma \times U_{\sigma'}$, where $T_\sigma
\cong \spec{k[\sigma^\bot_M]}$ is the minimal orbit of $U_\sigma$, and $U_{\sigma'}$
is the affine toric variety associated to the subsemigroup $\sigma_M' = \sigma_M /
\sigma^\bot_M$ of $M / \sigma^\bot_M$. Below, every construction with respect to
$(M, \leq_\sigma)$ will up to natural equivalence only depend on
$M_{\lessgtr_\sigma}$, and so for clearer presentation we will always neglect the
factor $\sigma_M^\bot$ and identify any $m \in M$ with its image in $\Z^{\sigma(1)}$.

The embedding of $M$ in $\Z^{\sigma(1)}$ is in a natural way compatible with the
partial order $\leq$ on $\Z^{\sigma(1)}$ induced by the subsemigroup
$\N^{\sigma(1)}$, i.e. $m \leq_\sigma m'$ iff $m \leq m'$. We consider the order
$\leq$ as an {\em extension} of $\leq_\sigma$ to $\Z^{\sigma(1)}$. For any $\tau <
\sigma$, the localization of $\leq_\sigma$ by $\tau_\sigma$ extends to a localization
of $\leq$ by the preorder $\leq'$ induced by the subsemigroup $\N^{\tau(1)} \oplus
\Z^{\sigma(1) \setminus \tau(1)}$, and we have a natural identification
$(\Z^{\sigma(1)})_{\lessgtr'} = \Z^{\tau(1)}$. This localization is naturally
compatible with the localization of $M$ by $\leq_\sigma^\tau$ and we have the
following commutative exact diagram:
\begin{equation*}
\xymatrix{
0 \ar[r] & \sigma^\bot_M \ar@{ >->}[d] \ar[r] & M \ar@{=}[d] \ar[r] & \Z^{\sigma(1)}
\ar@{-{>>}}[d]^{\pi} \ar[r] & A^\sigma \ar[r] \ar[d] & 0 \\
0 \ar[r] & \tau^\bot_M \ar[r] & M \ar[r] & \Z^{\tau(1)} \ar[r] & A^\tau \ar[r] & 0,
}
\end{equation*}
where $\pi$ is the canonical projection from $\Z^{\sigma(1)}$ onto $\Z^{\tau(1)}$.
Having these natural compatibilities in mind, in the sequel we will use the notation
$\leq_\sigma$ for both preorders on $M$ and on $\Z^{\sigma(1)}$; we will write
$\n \leq_\sigma m$ and the like for $\n \in \Z^{\sigma(1)}$ and $m \in M$.

We now describe more precisely the relation between the $\Delta$-family of \msh{E}
and the
$\hat{\Delta}$-family of $\pi^*\sh{E}$. For any $\sigma \in \Delta$ consider the
quotient representation $\pi_\sigma : U_{\hat{\sigma}} \twoheadrightarrow U_\sigma$,
where $U_{\hat{\sigma}} \cong k^{\sigma(1)}$. Here without loss of generality we
assume for the moment that $\sigma$ has full dimension in $N_\R$.
Denote $E^\sigma := \Gamma(U_\sigma, \sh{E})$ and $\hat{E}^\sigma := \Gamma(
U_{\hat{\sigma}}, \pi_\sigma^* \sh{E})$. The homogeneous coordinate ring has a
natural $A^\sigma$-grading $S^\sigma = \bigoplus_{\alpha \in A^\sigma}
S_\alpha^\sigma$, with $S_0 \cong k[\sigma_M]$, so that we can write:
\begin{equation*}
\hat{E}^\sigma \cong E^\sigma \otimes_{S^\sigma_0} S^\sigma \cong E^\sigma
\otimes_{S^\sigma_0} \big(\bigoplus_{\alpha \in A^\sigma} S_\alpha^\sigma\big) \cong
\bigoplus_{\alpha \in A^\sigma} \big(E^\sigma \otimes_{S^\sigma_0} S^\sigma_\alpha
\big).
\end{equation*}
By $E^\sigma \otimes_{S^\sigma_0} S^\sigma_0 \cong E^\sigma$, we find that
$E^\sigma$ is naturally embedded in $\hat{E}^\sigma$, and thus the $\sigma$-family
of $E^\sigma$ is a subfamily of the $\hat{\sigma}$-family of $\hat{E}^\sigma$.

Denote $\leq_\Delta$ the preorder on \mweildivisors, then every $(\Z^{\sigma(1)},
\leq_\sigma)$ is isomorphic to the localization of $(\Z^\rays, \leq_\Delta)$ by the
preorder $\leq'_\sigma$ where $\n \lessgtr'_\sigma \n'$ iff $\n - \n' \in
\Z^{\rays \setminus \sigma(1)}$. By the natural projection $\weildivisors
\longrightarrow \Z^{\sigma(1)}$, every fine graded module over the localization
$S_{x^{\hat{\sigma}}} \cong k[\Z^{\rays \setminus \sigma(1)} \times \N^{\sigma(1)}]$
is equivalent to a fine graded module over $S^\sigma \cong k[\N^{\sigma(1)}]$. The
localization
$\hat{E}_{x^{\hat{\sigma}}}$ of $\hat{E}$ then is equivalent to a representation
of $(\Z^{\sigma(1)}, \leq_\sigma)$, and by naturality,
the $\Delta$-family $E^\Delta$ glues as a subfamily of the $\hat{\Delta}$-family
$\hat{E}^{\hat{\Delta}}$.

\paragraph{Resolutions of $\hat{E}$.}
As for every equivariant coherent sheaf \msh{E} we can consider its associated fine
graded module $\hat{E}$, in principle there is nothing which prevents us from doing
this and to compute some finite free resolution of $\hat{E}$ over $S$, which then
descends to a resolution of \msh{E} of the desired type:
\begin{equation*}
0 \longrightarrow \breve{F}_s \longrightarrow \cdots \longrightarrow \breve{F}_0
\longrightarrow \sh{E} \longrightarrow 0,
\end{equation*}
where $\breve{F}_i \cong \bigoplus_{j = i}^{k_j} \sh{O}_X(D_{\n_{ij}})$ and the
length $s$ by Hilbert's syzygy theorem being bounded by the numbers of rays in
$\Delta$. At this point we could stop with this paper and leave the problem as an
application of traditional methods. However, there are some drawbacks of this
point of view, which motivate our further investigations. One problem is that the
pullback of a coherent sheaves along good quotients so far seems not to be
understood very well, not even for the toric case --- and as we will see in example
\ref{pullbacktorsionexample} below, such pullbacks might show pathological behaviour,
such as acquiring additional torsion. Another problem is that due to its global
nature, the module $\hat{\sh{E}}$ contains much more relations which might be
irrelevant to consider for getting a resolution.

\begin{example}
\label{pullbacktorsionexample}
Consider the subsemigroup $\sigma_M$ of $M \cong \Z^2$ which is generated by the
elements $(1, 0)$, $(1, 1)$, $(1, 2)$ and its associated semigroup ring $S_0 =
k[\sigma_M]$. Its fan is spanned by the primitive vectors $\n_1 = (2, -1)$ and
$\n_2 = (0, 1)$ in
$N_\R$, and the homogeneous coordinate ring $S = S^\sigma$ is $\Z_2$-graded.
Denote $\n := (-1, 0)$ and consider the reflexive $S_0$-module $S_{(\n)} \cong S_1$.
For the pullback we have:
\begin{equation*}
S_{(\n)} \otimes_{S_0} (S_0 \otimes S_1) \cong (S_1 \otimes_{S_0} S_0) \oplus
(S_1 \otimes_{S_0} S_1) \cong S_1 \otimes (S_1 \otimes_{S_0} S_1).
\end{equation*}
To compute $(S_1 \otimes_{S_0} S_1)$, we directly evaluate it as an $M$-graded tensor
product. The module $S_{(\n)}$ is a $M$-graded, where
\begin{equation*}
S_{(\n), m} =
\begin{cases}
k & \text{ if } \langle m, \n_i \rangle \geq 0 \text{ for } i = 1, 2 \\
0 & \text{ else}.
\end{cases}
\end{equation*}
In degree $m$, $S_{(\n)} \otimes_{k[\sigma_M]} S_{(\n)}$ is generated by all elements
$\chi(m_1) \otimes \chi(m_2)$ such that $m_1 + m_2 = m$ modulo the relation that
$\chi(m_1) \otimes \chi(m_2)$ is equivalent to $\chi(m_1 - m') \otimes \chi(m_2)$
and $\chi(m_1) \otimes \chi(m_2 - m'')$, respectively, whenever there exist some
$\chi(m')$ or $\chi(m'')$ such that $\chi(m_1 - m')$ and $\chi(m_2 - m'')$,
respectively, are in $S_{(\n)}$. It turns out that these relations cancel most of
the generators in every degree, so that for all nonzero degrees:
\begin{equation*}
\dim(S_{(\n)} \otimes_{S_0} S_{(\n)})_m =
\begin{cases}
2 & \text{ if } m = (1, 1) \\
1 & \text{ else}.
\end{cases}
\end{equation*}
Note that the nonzero degrees are precisely thos contained in the intersection of the
half spaces $\langle m, \n_2 \rangle
\geq 0$, $\langle m, 2\n_1 \rangle \geq 0$ and $\langle m, (1, 0) \rangle \geq 0$.
So, in degree $(1,1)$, our module has dimension two, whereas in all other degrees
it has at dimension one, which implies that it has torsion in degree $(1, 1)$, as
for any character $(1, 1) \leq_\sigma m$ the homomorphism $\chi_{(1, 1), m}$ can not
be injective. This is indeed an example where pullback of a torsion free, and even
reflexive, module along a geometric quotient aqcuires some new torsion.
\end{example}

Another phenomenon which we want to mention is that there are also other relevant
effects which one has to consider if one tries to choose some alternative module
instead of $\hat{E}$ whose descend coincides with that of $\hat{E}$.
For instance, for any affine toric variety $U_\sigma$ for which $A^\sigma$ is
nontrivial, there exist nonzero $S$-modules $F$ whose zero component vanishes;
the most easiest example is the one-dimensional module $S^\sigma / \langle x_\rho \mid
\rho \in \sigma(1) \rangle$ whose degree gets shifted by some nonzero $\alpha \in
A^\sigma$.

\section{Compression and Resolutions}
\label{resolutions}

\subsection{$\lcm$-lattices in $\Z^r$}

The partial order on $\Z^r$ induced by $\N^r$ coincides with the partial order given
by componentwise ordering, i.e. if we write $\n = (n_1, \dots, n_r)$, $\n' = (n'_1,
\dots, n'_r)$, then $\n \leq \n'$ iff $n_i \leq n'_i$ for every $1 \leq i \leq r$.
We set $\bar{\Z} := \{-\infty\} \cup \Z$ which is totally ordered by $-\infty < n$
for all $\n \in \Z$. Like $\Z^r$, the set $\bar{\Z}^r$ is partially ordered by the
componentwise total order, and the the canonical inclusion $\Z^r \hookrightarrow
\bar{\Z}^r$ is order preserving. We call any element
in $\bar{\Z}^r \setminus \Z^r$ {\em infinitary}.

For any element $\n \in \Z^r$ we can consider the subset $\n + \N^r$, which is
the intersection of the shifted cone $\n + \R_{\geq 0}^r$ with $\Z^r$. It is easy to
see that for any finite set of elements $\n_1, \dots \n_s$ in $\Z^r$, the intersection
$\bigcap_{i = 1}^s \big(\n_i + \N^r\big)$ is again of the form $\n + \N^r$. The
element $\n$ is called the {\em least common multiple} of $\n_1, \dots, \n_r$, denoted
$\lcm\{\n_1, \dots, \n_s\}$, and it is given by componentwise maximum of the $\n_i$.
The $\lcm$ extends canonically to $\bar{\Z}^r$. In the geometric picture, for some
infinitary element $\n = (n_1, \dots, n_r)$ with $n_{i_j} = -\infty$ for some
$\{i_1, \dots, i_r\} \subset \{1, \dots, r\}$, we write $\n + C$ for the cone, where
$C = \{ \underline{c} \in \R^r \mid c_i \geq 0 \text{ if } i \notin \{i_1, \dots,
i_k\}\}$. One can think of the cone $C$ of the standard orthant moved to minus
infinity in the directions $i_1, \dots, i_k$.

In our actual definition of the $\lcm$-lattice, we will need inifinitary elements
to generate the lattice, but after generation, we throw away all these elements.
Instead, we close every $\lcm$-lattice from below by adding the unique minimal
element $(-\infty, \dots, -\infty) =: \hat{0}$.

\begin{definition}
Let $\p \subset \bar{\Z}^r$ be some poset and $\lcm(\p)$ the lattice generated by
the $\lcm$'s of elements in $\p$. Then we denote the set $(\lcm(\p) \cap \Z^r)
\cup \hat{0}$ the {\em $\lcm$-lattice} of $\p$.
\end{definition}

Every $\lcm$-lattice $\cL$ gives rise to a partition of $\Z^r$, respectively to an
equivalence relation,
on $\Z^r$. Namely, for every element $\underline{n} \in \Z^r$, there exists a unique
maximal element $\n' \in \cL$ with $\n' \leq \n$.

\begin{definition}
Let $\n \in \Z^r$ and $\n' \in \cL$ maximal such that $\n' \leq \n$. Then we call
$\n'$ the {\em anchor element} $A(\n)$ of $\n$ in $\cL$. Any two elements $\n_1, \n_2
\in \Z^r$ are equivalent iff $A(\n_1) = A(\n_2)$. We denote $T_\n$ the equivalence
class associated to $\n \in \cL$.
\end{definition}

\subsection{Polynomial rings}
\label{polynomialrings}

In this subsection we consider the special case where $X$ is an affine toric variety
isomorphic to the affine space $k^r$, so that we can assume without loss of
generality that $\sigma$ and $\sigma_M$ coincide with the standard orthant
$\R^r_{\geq 0}$ in $\R^r$, and the subsemigroup $\N^r$ of $\Z^r$, respectively.
We denote $S \cong k[\N^r]$ the coordinate ring of $X$ and $E$ a nonzero finitely
generated $S$-module. We formally extend the representation of  $(\Z^r, \leq)$ by $E$
to a representation of $(\bar{\Z}^r, \leq)$ by setting $E_\n = 0$ for all infinitary
$\n$. In order
to construct a compression functor for $E$, we have to extract all nontrivial maps
(i.e. the nonisomorphisms) of the corresponding $\sigma$-family, as well as all
possible relations among them.

\begin{definition}
Let $\n \in \Z^r$, then we define the set $I_E(\n)$ to contain those elements $\n'$
in $\bar{\Z}^r$ which are minimal with the property that for all $\n'' \in
\bar{Z}^r$ with $\n' \leq \n'' \leq \n$ the morphisms $\chi_{\n', \n''}: E_{\n'}
\rightarrow E_{\n''}$ and $\chi_{\n'', \n}: E_{\n''} \rightarrow E_{\n}$ are
isomorphisms. We denote $\mathcal{I}_E := \bigcup_{\n \in \Z^r} I_E(\n)$.
\end{definition}

Note that the case where $I_E(\n)$ contains an infinitary elemement can only
(but not necessarily has to) occur when $E_\n$ is zero. Moreover, note that it
follows immediately from the finitely generatedness of $E$ that
$\mathcal{I}_E$ and the $I_E(\n)$ are finite sets.

\begin{definition}
We denote $\cL_E$ the $\lcm$-lattice generated by $\mathcal{I}_E$.
For any $\n \in \Z^r$, we denote the corresponding anchor element by $A_E(\n)$.
\end{definition}

We can depict the set of equivalence classes as a tiling of $\R^r$ by cubic, possibly
non-compact blocks, where the anchor elements are precisely those elements sitting on
the smallest vertex with respect to $\leq$.
Observe that $\lcm\{\n_1, \dots, \n_s\} \in \Z^r$ as soon as at least one of the
$\n_i$ is non-infinitary. Moreover, if $I_E(\n)$ contains an infinitary element, this
implies that $E_\n = 0$. In general, the set $I_E(\n)$ will contain
infinitary elements only if there exists no $\n' < \n$ such that $E_{\n'} \neq 0$.
In that case, $I_E(\n)$ will contain $\hat{0}$ as its only element. An exception
are those modules $E$, which are of rank zero, and thus are torsion modules.
The infinitary elements $I_E(\n)$ for all $\n \in \Z^r$ in that case describe the
support of $E$.

\begin{example}
Let $J \subset S$ be a monomial ideal, generated by monomials $x^{\n_1}, \dots,
x^{\n_s}$. Then we have
\begin{equation*}
I_J(\n) =
\begin{cases}
\{\n_i \leq \n \} & \text{ if } x^\n \in J \\
\hat{0} & \text{ else},
\end{cases}
\end{equation*}
and the anchor element $A_J(\n)$ being $\lcm\{\n_i \leq \n\}$. The lattice $\cL_J$
then coincides with the $\lcm$-lattice introduced in \cite{GPW99}.
\end{example}

\begin{figure}[hbpt]
\centerline{\input{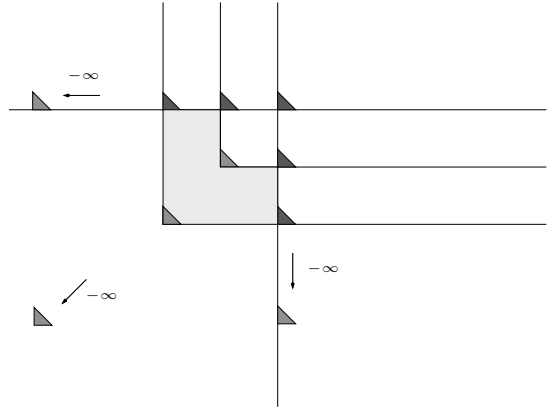}}
\caption{$\lcm$-lattice for example \ref{lcmexample1}\label{f-lcmexample1}}
\end{figure}

\begin{example}
\label{lcmexample1}
Consider the torsion module $T= k[x, y] / \langle x^2, xy, y^2\rangle$. We have
\begin{equation*}
I_T(\n) =
\begin{cases}
\{(0, 0)\} & \text{ for } \n \in \{(0, 0), (1, 0), (0, 1)\} \\
\{(1, 1)\} & \text{ for } \n = (1, 1) \\
\{(2, -\infty)\} & \text{ for } \n = (k, 0), k > 1 \\
\{(-\infty, 2)\} & \text{ for } \n = (0, k), k > 1 \\
\{(1, 1), (2, -\infty)\} & \text{ for } \n = (k, 1), k > 1 \\
\{(1, 1), (-\infty, 2)\} & \text{ for } \n = (1, k), k > 1 \\
\{(1, 1), (2, -\infty), (-\infty, 2)\} & \text{ for } (2, 2) \leq \n \\
\{\hat{0}\} & \text{ else}.
\end{cases}
\end{equation*}
The corresponding $\lcm$-lattice then is the set $\{\hat{0}, (0, 0), (2, 0),
(0, 2), (1, 1), (1, 2),$ $(2, 1),$ $(2, 2)\}$. Figure \ref{f-lcmexample1} shows
the partitioning of $\Z^2$ by the $\lcm$-lattice. The rectangular figure indicates
the degrees $(0, 0), (1, 0), (0, 1)$, where $T$ is nonzero; the light grey
triangles indicate all the initial elements $I_T(\n)$, and the darker grey triangles
denote the additional elements of the $\lcm$-lattice. The infinitary
elements become merged to $\hat{0}$ in $\cL_T$.
\end{example}

Denote $\cL_E$-Rep the category of finite-dimensional $k$-linear representations
of $\cL_E$; denote $\mathcal{M}_E$ the full subcategory of the category of
fine-graded $S$-modules whose objects are those modules $F$ whose associated
$\lcm$-lattice
$\cL_F$ is a sublattice of $\cL_E$. Let $\iota_E : \cL_E \hookrightarrow \Z^r$ be
the canonical inclusion. Then we define the functor $\zip^E$ from $\mathcal{M}_E$
into $\cL_E$-Rep by
\begin{equation*}
\zip^E(F) := \iota_E^* F
\end{equation*}
where $\iota_E^*$ denotes the sheaf pullback.

To define the $\unzip$ functor, we have to do a little bit more. Let $F$ be some
representation of $\cL_E$, mapping $\n$ to $F(\n)$, and $\n \leq \n'$ to $F(\n, \n')$.
Then we define a representation of $\Z^r$ by setting $F_\n := F\big(A_E(\n)\big)$ and
$\chi_{\n, \n'} := F\big(A_E(\n), A_E(\n')\big)$ for every pair $\n, \n' \in \Z^r$.
This indeed establishes a well defined functor, where $F\big(A_E(\n),
A_E(\n')\big) = \operatorname{id}$ whenever $A_E(\n) = A_E(\n')$ and $F\big(A_E(\n),
A_E(\n'')\big) = F\big(A_E(\n'), A_E(\n'')\big) \circ F\big(A_E(\n), A_E(\n')\big)$
whenever $\n \leq \n' \leq \n''$.

\begin{theorem}
The pair of functors $\zip$ and $\unzip$ establishes an equivalence of categories
between $\mathcal{M}_E$ and $\cL_E$-Rep.
\end{theorem}
 
\begin{proof}
We show that $\unzip \circ \zip \cong 1_{\mathcal{M}_E}$ and $\zip \circ \unzip
\cong 1_{\cL_E\operatorname{-Rep}}$. In the first case, let $F$ be some
representation of $(\Z^r, \leq)$. Denote $F'_\n := \unzip(\iota_E^*F)(\n)$ for
every $\n \in \Z^r$ and define $h : F'_\n \longrightarrow F_\n$ by setting
$h := \chi_{A_E(\n), \n}$. Now $h$ is an isomorphism for every $\n \in \Z^r$, and
moreover, for every pair $\n \leq \n'$, we have $\chi_{A_E(\n'), \n'} \circ
\chi_{A_E(\n), A_E(\n')} = \chi_{\n, \n'} \circ \chi_{A_E(\n), \n} = \chi_{A_E(\n),
\n'}$. So we obtain $\unzip \circ \zip \cong 1_{\mathcal{M}_E}$.

The other direction is immediate, and we even obtain $\zip \circ \unzip =
1_{\cL_E\operatorname{-Rep}}$
\end{proof}

\begin{corollary}
$\mathcal{M}_E$ is an abelian category.
\end{corollary}

Let $\n$ be any element in $\cL_E$, then we can consider the free representation
$F^\n$ of $\cL_E$. Its unzipping has a particularly easy structure, namely
$\unzip(F^\n) \cong S(-\n)$, i.e. the free fine-graded $S$-module with degree
shifted by $-\n$. $\unzip(F^\n)$ is the unique $S$-module which has the property
that its $\n'$-th degree is one-dimensional if $\n \leq \n'$ and zero else.

Now we can consider a free resolution of $\zip(E)$ in terms of free representations
of $\cL_E$:
\begin{equation*}
0 \longrightarrow F_s \longrightarrow \cdots \longrightarrow F_0 \longrightarrow
\zip(E) \longrightarrow 0
\end{equation*}
where for every $1 \leq i \leq s$:
\begin{equation*}
F_i \cong \bigoplus_{\n \in \cL_E} (F^\n)^{f^i_\n}
\end{equation*}
where $f^i_\n$ is the free dimension of the vector space associated to $\n$ in the
$(i - 1)$-th syzygy representation.
By unzipping, we obtain an exact sequence of fine-graded $S$-modules:
\begin{equation}
\label{Sres}
0 \longrightarrow \unzip(F_s) \longrightarrow \cdots \longrightarrow \unzip(F_0)
\longrightarrow E \longrightarrow 0
\end{equation}
where for every $1 \leq i \leq s$:
\begin{equation*}
\unzip(F_i) \cong \bigoplus_{\n \in \cL_E} S(-\n)^{f^i_\n}.
\end{equation*}
In order to show, that this is a minimal free resolution of $E$ over $S$, we consider
the first step of the resolution $0 \rightarrow K_0 \rightarrow \unzip F_0
\rightarrow E
\rightarrow 0$. We define a map $\phi : \cL_E \longrightarrow \cL_{K_0}$ by mapping
every $\n \in \cL_E$ to its anchor element in $\cL_{K_0}$:
\begin{equation*}
\phi(\n) := A_{K_0}(\n).
\end{equation*}
We have the following:
\begin{proposition}
\label{syzcontraction}
The map $\phi$ is a contraction.
\end{proposition}

\begin{proof}
We first show that $\phi(U_E(\n)) = U_{K_0}(A_{K_0}(\n))$ for all $\n \in \cL_E$,
where we write $U_E$ and $U_{K_0}$ for open subsets in $\cL_E$ and $\cL_{K_0}$,
respectively. Clearly, $\phi(U_E(\n)) \subset U(A_{K_0}(n))$; by construction of
$K_0$, the lattice $\cL_{K_0}$ is a sublattice of $\cL_E$, so that for any $\n' \in
U_{K_0}(A_{K_0})$ there is $\n'' \in U_E(\n)$ with $\phi(\n'') = \n'$.
Now let $\n \in \cL_{K_0}$ and consider the set $\phi^{-1}\big(U_{K_0}(\n)\big)$,
which consists of all $\n' \in \cL_E$ such that $\n \leq A_{K_0}(\n')$. $\n \leq \n'$
implies $\n = A_{K_0}(\n) \leq A_{K_0}(\n')$, and thus $U_E(\n) \subset
\phi^{-1}(U_{K_0}(\n))$. Moreover, $\phi^{-1}\big(U_{K_0}(\n)\big) = \{\n' \in \cL_E
\mid \n \leq A_{K_0}(\n')\}$, and thus $\phi^{-1}\big(U_{K_0}(\n)\big) \subset
U_E(\n)$. Hence, $\phi^{-1}\big(U_{K_0}(\n)\big) = U_E(\n)$, and $\phi$ is a
contraction.
\end{proof}

\begin{theorem}
Sequence (\ref{Sres}) is a minimal free resolution of $E$ over $S$.
\end{theorem}

\begin{proof}
Observe that the number of $k$-linear independent generators of the module $E$
degree $\n$ is the codimension of the subvector space $\sum_{\n' < \n} x^{\n - \n'}
\cdot E_{\n'}$ of $E_\n$, which coincides with the free dimension of $E_\n$. Thus
$\unzip F_0$ is the minimal free module which surjects onto $E$.
Using proposition \ref{syzcontraction} and lemma \ref{contractionliftres}, we see
that a resolution of $K_0$ over $\cL_E$ is a lift of some resolution of $K_0$
restricted to $\cL_{K_0}$. Hence, the theorem follows by induction.
\end{proof}

\subsection{Admissible posets and normal semigroup rings}
\label{semigrouprings}

To extend  our considerations to the case of normal semigroup rings, consider the
map $M \rightarrow \Z^{\sigma(1)}$,
which without loss of generality we assume to be injective. This corresponds to a
quotient representation $\pi : k^{\sigma(1)} \twoheadrightarrow U_\sigma$ together
with an $A$-graded homogeneous coordinate ring $S := k[x_\rho \mid \rho \in
\sigma(1)]$. For any coherent sheaf $\sh{E}$ over $U_\sigma$, we can consider its
pullback $\pi^* \sh{E}$ over $k^{\sigma(1)}$. 

Applying the machinery from subsection \ref{polynomialrings}, we can obtain a
reflexive resolution for $\sh{E}$ by sheafification of the resolution of $\hat{E}$
with respect to the $\lcm$-lattice $\cL_{\hat{E}}$:
\begin{equation*}
0 \longrightarrow \unzip(F_r)\breve{\ } \longrightarrow \cdots \longrightarrow
\unzip(F_0)\breve{\ } \longrightarrow \sh{E} \longrightarrow 0,
\end{equation*}
where $\sh{E} \cong \big(\unzip \iota_{\hat{E}}^* \pi^* E\big)\breve{\ }$.
For any anchor element $\n \in \cL_{\hat{E}}$, the unzipping of the associated free
representation of $\cL_{\hat{E}}$ is isomorphic to $S(-\n)$.
Unlike the case of smooth toric varieties, in the general case such a
resolution is not uniquely defined, and it can be possible to obtain shorter
resolutions which are of this type.

\begin{definition}
\label{admissibledef}
Let $E$ be a $M$-graded $k[\sigma_M]$-module. A finite subposet $\p \subset \Z^r \cup
\hat{0}$ is {\em admissible with respect to $E$} if
\begin{enumerate}[(i)]
\item\label{admissibledefi} for all $m \in M$ there exists a {\em unique} $\n \in
\p$ with $\n \leq m$, such that $\n' \leq m$ implies $\n' \leq \n$ for all $\n'
\in \p$;
\item\label{admissibledefii} consider the open set $U_\n = \bigcup_{\n \leq m} U(m)$
in $M$ and the vector space $E(U_\n) = \underset{\leftarrow}{\lim} E_m$, there exists
a vector space $E_\n$ and a diagonal homomorphism $E_\n \longrightarrow E(U_\n)$ such
that every induced homomorphism $E_\n \longrightarrow E_m$ is an isomorphism for all
$m \in T_\n$.
\end{enumerate}
We call $E_\n$ the {\em anchor completion} of $E$ at $\n$ and we denote $A_E(m)$ the
unique maximal element $\n \in \p$ with $\n \leq m$.
\end{definition}
Note that in the definition we have identified the elements $m \in M$ with their
image in $\Z^{\sigma(1)}$. For any $\n \in \p$, the homomorphism $E_\n \rightarrow
E(U_\n)$ necessarily is injective, and for every $\n \leq \n'$, the composition
\begin{equation*}
E_\n \longrightarrow E(U_\n) \longrightarrow E(U_{\n'})
\end{equation*}
is a diagonal morphism, which factors through the image of $E_{\n'}$, such that we
obtain a morphism between the anchor completions $E_\n \longrightarrow E_{\n'}$.

\begin{lemma}
Assume that $U_\sigma$ is smooth and thus $M \cong \Z^{\sigma(1)}$ and let $\p$ be
some admissible poset with respect to $E$. Then for every subset $m_1, \dots, m_s$ of
$\p$, $\lcm\{m_1, \dots, m_s\}$ is also contained in $\p$. In particular, $\p$
contains the $\lcm$-lattice $\cL_E$.
\end{lemma}

\begin{proof}
Denote $m_l := \lcm\{m_1, \dots, m_s\}$. There exists a unique $m \in \p$ such that
$m \geq_\sigma m_l$; but such an $m$ must coincide with $m_l$.
\end{proof}

\comment{
In the case where $\sigma$ has not full dimension, we can consider a splitting of
$M \cong M^\bot_\sigma \oplus M_\sigma$, where $M_\sigma$ is the image of $M$ in
$\Z^{\sigma(1)}$ by the map $M \rightarrow \Z^{\sigma(1)}$ and $M_\sigma^\bot$ its
kernel. $(M_\sigma) \otimes_\Z \R$ the can naturally be identified with the dual
space of the minimal subspace $N_\sigma \subset N_\R$ containing $\sigma$. Moreover,
as partially ordered set, $M_\sigma$ is isomorphic to $M_{\lessgtr_\sigma}$, and
any representation of $M_\sigma$ with respect to the induced partial order is
equivalent to a representation of the preordered set $(M, \leq_\sigma)$.
}

From the observation that $\cL_{\hat{E}}$ is admissible, we conclude:

\begin{proposition}
Every finitely generated $k[\sigma_M]$-module $E$ has an admissible poset.
\end{proposition}

\begin{proof}
We take the poset of all $\n \in \cL_{\hat{E}}$ such that $\{m \in M \mid A_E(m)
= \n\} \neq \emptyset$.
\end{proof}

Let $\p \subset \Z^{\sigma(1)}$ be an admissible poset, and denote $\mathcal{M}_\p$
the category of finitely generated, $M$-graded $k[\sigma_M]$-modules for which $\p$
is admissible.
Then we define the functor $\zip^\p$ from $\mathcal{M}_\p$ to the category of
$k$-linear representations of $\p$ by:
\begin{equation*}
\zip^\p (E)_\n := E_\n,
\end{equation*}
where $E_\n$ is the anchor completion at $\n$.

\begin{remark}
\label{bigadmissible}
Our definition also allows to add anchor elements $\n$ such that the corresponding
set $T_\n$ is empty. In that case we set $E_\n = \underset{\leftarrow}{\lim} E_{\n'}$
for all $\n < \n' \in \p$ such that $T_{\n'} \neq \emptyset$.
\end{remark}

In the opposite direction, from every representation $E$ of an admissible poset $\p$
one can construct a representation of $M$. We define $\unzip^\p (E)$ by setting:
\begin{enumerate}[(i)]
\item $\unzip^\p (E)_m := E_{A(m)}$,
\item $\chi_{m, m'} := E\big(A(m), A(m')\big)$.
\end{enumerate}

\begin{theorem}
\label{admissibleequivalence}
The pair $\zip^\p$ and $\unzip^\p$ is a compression of $\mathcal{M}_\p$, i.e.
$\zip^\p$ and $\unzip^\p$ are functors which establish an equivalence of categories.
\end{theorem}

\begin{proof}
By construction, $\unzip^\p \circ \zip^\p (E) \cong E$ for every $k[\sigma_M]$-module
for which $\p$ is admissible. To obtain functors, we show that any morphism $E
\longrightarrow F$ of objects in $\mathcal{M}_\p$ induces a morphism of the
corresponding representations of $\p$ and vice versa. First, any homomorphism $E
\rightarrow F$ is a homomorphism of sheaves over $(M, \leq_\sigma)$, and thus there
is an induced homomorphism $E_\n \rightarrow E(U_\n) \rightarrow F(U_\n)$ for every
$\n \in \p$, which factors through the diagonal $F_\n$, hence we obtain a
homomorphism $E_\n \rightarrow F_\n$; the family of such morphisms for every $\n \in
\p$ in a natural way represents a homomorphism of representations of $\p$. In the
other direction, a homomorphism $f: \zip^\p (E) \rightarrow \zip^\p(F)$ unzips
componentwise as $f_m := f_{A(m)} : E_{A(m)} \rightarrow F_{A(m)}$.
\end{proof}

\begin{proposition}
Let $\n \in \p$, then $\p$ is admissible with respect to the reflexive module
$S_{(\n)}$, and moreover, $S_{\n} \cong \unzip^\p F^\n$.
\end{proposition}

\begin{proof}
Let $m \in M$, then $\n \leq A_E(m)$ iff $\n \leq m$: the first implication is clear,
because $m \leq A_E(m)$; for the second, observe that $A_E(m) \geq \lcm \{A_E(m),
\n\}$, and thus $\n \leq A_E(m)$. So $\p$ is admissible with respect to $S_{(\n)}$
and $\unzip^\p F^\n \cong S_{(\n)}$.
\end{proof}

As in the case for polynomial rings, we obtain a reflexive resolution for $E$:
\begin{equation*}
0 \longrightarrow \unzip^\p(F_s) \longrightarrow \cdots \longrightarrow
\unzip^\p(F_0) \longrightarrow E \longrightarrow 0.
\end{equation*}

\begin{figure}[ht]
\includegraphics[height=5cm,width=5cm]{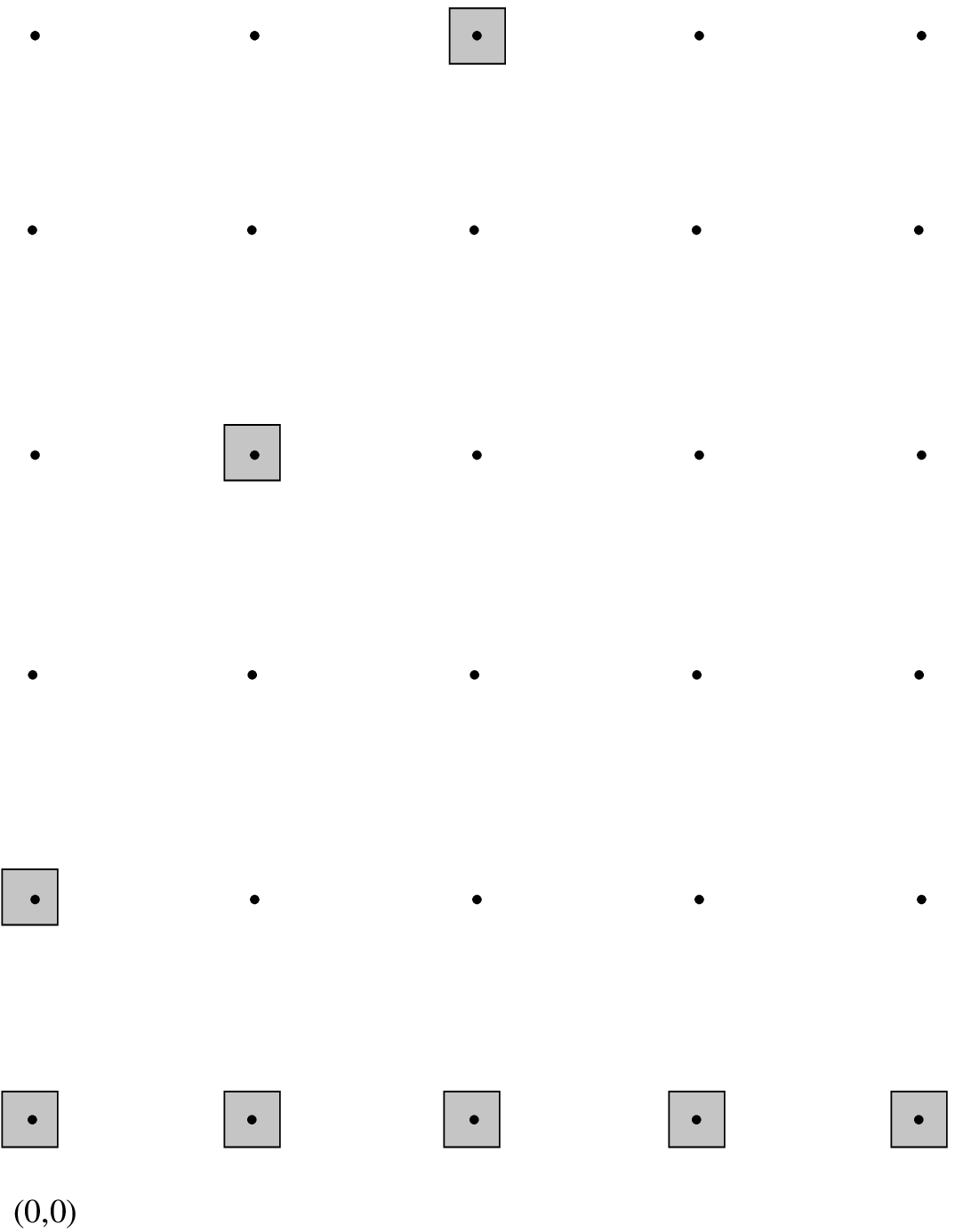}\quad\quad\quad\quad
\includegraphics[height=5cm, width=8cm]{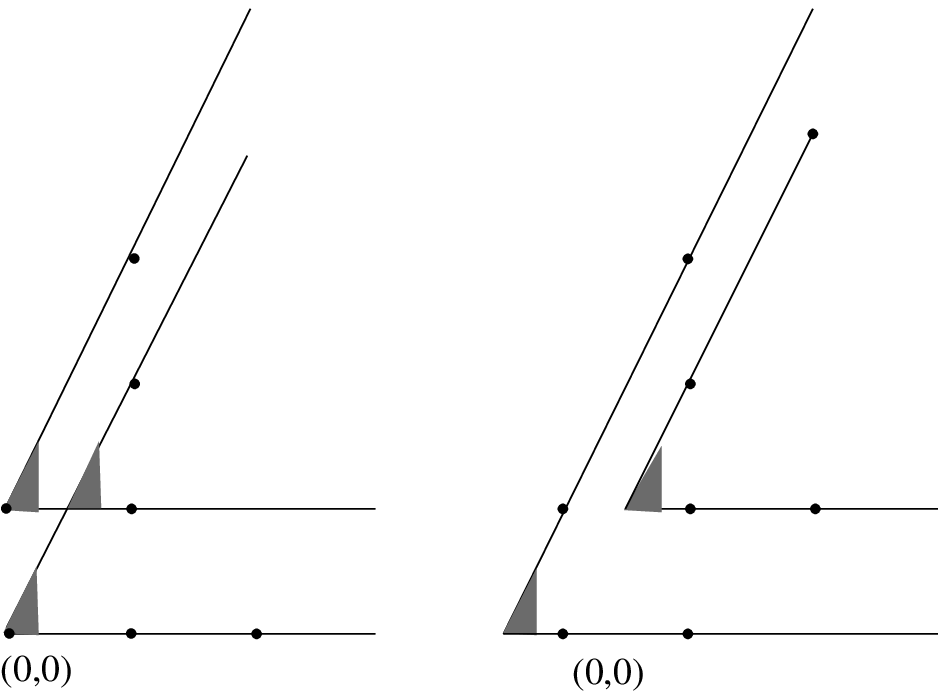}
\caption{The module from example \ref{admissibleexample1} and the partitions of
$\Z^2$ with respect to $\p_1$ and $\p_2$}\label{f-admissibleexample1}
\end{figure}

\begin{example}
\label{admissibleexample1}
Consider the semigroup $\sigma_M$ from \ref{pullbacktorsionexample} and the torsion
sheaf $T$ which is given by:
\begin{equation*}
T_m =
\begin{cases}
k & m = (p, 0), p \geq 0\\
k & m = (0, 1) + p \cdot (1, 2), p \geq 0\\
0 & \text{ else},
\end{cases}
\end{equation*}
and $\chi_{m, m'} = \operatorname{id}$ whenever $T_m, T_{m'} \neq 0$.
We can compare the following two admissible posets,
\begin{equation*}
\p_1 = \{\hat{0}, (0, 0), (-1, 1), (0, 1)\}
\end{equation*}
and
\begin{equation*}
\p_2 = \{\hat{0}, (-1, 0), (0, 1)\}.
\end{equation*}
We have $(\zip^{\p_1}T)_{(0, 0)} = k$, $(\zip^{\p_1}T)_{(-1, 1)} = k$,
$(\zip^{\p_1}T)_{(0, 1)} = 0$, and $\zip^{\p_2}T_{(-1, 0)} = k$,
$\zip^{\p_2}T_{(0, 1)} = 0$. But in the latter case, we have that $T(U_{(-1, 0)})$ is
is the fiber product $k \times_0 k \cong k^2$, such that $\zip^{\p_2}T_{(-1, 0)}$
corresponds to a proper
diagonal homomorphism $k \rightarrow k^2$. These compressions give rise to two
somewhat different resolutions. Via resolving over $\p_1$ and by $\unzip^{\p_1}$, we
obtain:
\begin{equation*}
0 \longrightarrow S_{(0, -1)} \longrightarrow S_{(1, -1)} \oplus S_{(0, 0)}
\longrightarrow T \longrightarrow 0
\end{equation*}
and for $\p_2$:
\begin{equation*}
0 \longrightarrow S_{(0, -1)} \longrightarrow S_{(1, 0)} \longrightarrow T
\longrightarrow 0.
\end{equation*}
In a sense, the module $T$ is like the module $k[x, y] / \langle xy \rangle$
over the polynomial ring $k[x, y]$, whose $\lcm$-lattice is isomorphic to $\p_2$.
However, there exists no unique minimal element in $m \in M$ with $T_m \neq 0$,
so that the consideration of the diagonal morphism indeed is necessary to obtain
a resolution which is like the minimal resolution of $k[x, y] / \langle xy \rangle$.
The left part of figure \ref{f-admissibleexample1} shows a part of the lattice
$\Z^2$; the light grey areas indicate the degrees, where $T_m$ is nonzero. The right
part of figure \ref{f-admissibleexample1} shows the partitioning of $\Z^2$ according
to the two admissible posets $\p_1$ and $\p_2$.
\end{example}

\subsection{Extension of a module to the homogeneous coordinate ring}
\label{homext}

Consider $E$ any $M$-graded $k[\sigma_M]$-module and $S$ the homogeneous coordinate
ring for $U_\sigma$. As we have seen in the previous subsection, one can in a natural
way associate the $S$-module $\hat{E} = E \otimes_{S_0} S$ to $E$. In this subsection
we want to discuss another way to associate a module, denoted $EE$, to $E$ which also
has the property that $E\breve{E} \cong E$, but which behaves better, for instance it
preserves
the property of torsion freeness. For this, for every $\n \in \Z^{\sigma(1)}$ we
denote $U_\n := \bigcup_{\n \leq m} U(m)$ an open subset of $(M, \leq_\sigma)$. To
see that the set $\{\n \leq m\}$ is always nonempty, just choose some $m \in
\sigma_M$ with $\langle m, n(\rho) \rangle > 0$ for every $\rho \in \sigma(1)$, then
for every $\n \in \Z^{\sigma(1)}$ we can choose an integer $c > 0$ such that $\n \leq
c \cdot m$. Thus for every $\n \in \Z^{\sigma(1)}$ the vector space $E(U_\n)$ exists,
and can be identified with $\underset{\leftarrow}{\lim} E_m$. If $E$ is finitely
generated, then the $E(U_\n)$ are finite dimensional.

\begin{definition}
We define a representation of $(\Z^{\sigma(1)}, \leq)$ by:
\begin{equation*}
EE_\n := E(U_\n).
\end{equation*}
\end{definition}

For every $\n \leq \n'$, the set $U_\n'$ is contained in $U_\n$, and thus we have
a functorial homomorphism $EE_\n \rightarrow EE_{\n'}$, and indeed we obtain a
well-defined representation of $\Z^{\sigma(1)}$.

\begin{proposition}
\label{EEprop}
$EE$ has the following properties:
\begin{enumerate}[(i)]
\item\label{EEpropi} $E\breve{E} = E$.
\item\label{EEpropii} If $E$ is finitely generated, then also $EE$ is finitely
generated.
\item\label{EEpropiii} If $E$ is torsion free, then also $EE$ is torsion free.
\end{enumerate}
\end{proposition}

\begin{proof}
(\ref{EEpropi}): By definition, if $\n = m \in M$, then $EE_m = E(U(m)) = E_m$, thus
$EE_0 = E$.
(\ref{EEpropii}): We apply the criteria of \cite{perling1}, \S 5.3. We have already
stated that the $EE_\n$ are finite dimensional. For all infinite chains $\cdots <
\n_i < \n_{i + 1} < \cdots$, we know that there exists an index $i_0$ such that the
$E_m$ vanish for $m \leq \n_{i_0}$, and thus the $EE_\n$ are zero. To see that there
are only finitely many $\n$ such that $\bigoplus_{\n' < \n} EE_{\n'} \rightarrow
EE_\n$ is not surjective, we choose some finite poset in $\Z^{\sigma(1)}$ which is
admissible with respect to $E$; using this, we find that there are only finitely many
isomorphism classes of vector spaces $EE_\n$.
(\ref{EEpropiii}): As the morphisms $\chi_{m, m'}$ are injective for every $m
\leq_\sigma m'$, the induced morphisms of the limits $EE_\n \rightarrow EE_{\n'}$
are also injective for every $\n \leq \n'$.
\end{proof}

Note that in general $EE$ is not just $\hat{E}$ modulo torsion. For instance, the
module $\check{E}$ of example \ref{pullbacktorsionexample} modulo torsion is not
reflexive, whereas $EE$ is reflexive (see subsection \ref{reflext}).

\begin{proposition}
The $\lcm$-lattice of $EE$ is admissible with respect to $E$.
\end{proposition}

\begin{proof}
Denote $\mathcal{L}$ the $\lcm$-lattice of $E$.
Let $m \in M$ and $\n \in \mathcal{L}$ its anchor element. By definition, the
map $EE_\n \longrightarrow E_m$ is an isomorphism, and thus $\mathcal{L}$ is
admissible.
\end{proof}

So we can use $EE$ as alternative module by which we can construct resolutions of
\msh{E}. In subsection \ref{reflext}, we will do a more explicit analysis of $EE$
for the case where \msh{E} is reflexive.

\subsection{Global resolutions for $\Delta$-families}
\label{deltaglobres}

Now let \msh{E} be an equivariant coherent sheaf over $X$ and $E^\Delta$ its
associated $\Delta$-family. To obtain global resolution of \msh{E}, we want to
extend the techniques considered in the previous two subsections. Denoting
$E^\sigma := \Gamma(U_\sigma, \sh{E})$, we assume that we have a family
$\mathfrak{P} = \{\p^\sigma \mid \sigma \in \Delta\}$ of posets and
compressions $\zip^{\p^\sigma}, \unzip^{\p^\sigma}$ with respect to these posets.
For nicer notation, we write $\zip^\sigma$
and $\unzip^\sigma$ instead of $\zip^{\p^\sigma}$ and $\unzip^{\p^\sigma}$.
For any $m \in M$ we write $A_E^\sigma(m)$ for the anchor element of $m$ in
$\p^\sigma$. We denote $l_{\sigma\tau}$ and $k_{\sigma\tau}$ the gluing maps
for the families $(M, \leq_\sigma)$ and $\p^\sigma$.

\begin{definition}
The collection $\mathfrak{P} = \{\p^\sigma \mid \sigma \in \Delta\}$ is called
{\em admissible} with respect to \msh{E} if it glues over $\Delta$ and for every
$\sigma \in \Delta$, the poset $\p^\sigma$ is admissible with respect to $E^\sigma$.
\end{definition}

\comment{
For any $\tau < \sigma$
 and the corresponding admissible posets $\p^\tau$ and
$\p^\sigma$, by abuse of notation we denote $\leq_\tau$ and $\leq_\sigma$ the
partial orders on $\p^\tau$ and $\p^\sigma$, respectively, which in a natural way
are compatible with the corresponding preorders on $M$. The localization
in a natural way is compatible the localization of $\leq_\sigma$ by
$\leq_\tau$ on $M$. This way, via the canonical projection of $\Z^{\sigma(1)}$ onto
$\Z^{\tau(1)}$ we identify the set
$(\p^\sigma)_{\lessgtr_\sigma^\tau}$ with its image in $\Z^{\tau(1)}$.
}

\paragraph{Compressions of $\Delta$-families.}

\begin{proposition}
Let $\mathfrak{P} = \{\p^\sigma \mid \sigma \in \Delta\}$ be a collection of posets
which is admissible with respect to \msh{E} and assume that we have a family of
sheaves $F^\sigma$ which glues over the collection $\p^\sigma$, then the family
$\unzip^\sigma F^\sigma$ is a $\Delta$-family.
\end{proposition}

\begin{proof}
We show that $l^*_{\sigma\tau} \unzip^\tau F^\tau \cong \unzip^\sigma
k_{\sigma\tau}^* F^\tau$ for every $\tau < \sigma$. This follows componentwise from
$(l^*_{\sigma\tau} \unzip^\tau F^\tau)_m = (\unzip^\tau F^\tau)_{l_{\sigma\tau}(m)}
= (\unzip^\tau F^\tau)_m = F^\tau_{A^\tau(m)}$ and
$(\unzip^\sigma k_{\sigma\tau}^* F^\tau)_m = (k^*_{\sigma\tau} F^\tau)_{A^\sigma(m)}
= F^\tau_{k_{\sigma\tau}(A^\sigma(m))} \cong F^\tau_{A^\tau(m)}$, where the last
isomorphism follows from the fact that $k_{\sigma\tau}$ is a contraction.
Denote $\Psi^{\sigma\tau} : k_{\sigma\tau}^* F^\tau \overset{\cong}{\longrightarrow}
F^\sigma$ the gluing maps over the family $\p^\sigma$, then we set $\Phi^{\sigma\tau}
:= \unzip^\sigma \Psi^{\sigma_\tau}$.
By the isomorphisms $l^*_{\sigma\tau} \unzip^\tau F^\tau \overset{\cong}{\rightarrow}
\unzip^\sigma k_{\sigma\tau}^* F^\tau$ for all $\tau < \sigma$ and the functoriality
of $l_{\sigma\tau}^*$, we have for any triple $\rho < \tau < \sigma$ the natural
identification $\Phi^{\sigma\rho} = \Phi^{\sigma \tau} \circ l_{\sigma\tau}^*
\Phi^{\tau\rho}$, and the proposition follows.
\end{proof}

Denote $\mathbf{S}^{\mathfrak{P}}$ the category of coherent equivariant sheaves over
$X$ with respect to which the collection $\p^\sigma$ is admissible.
The operations $\zip^\Delta$ and
$\unzip^\Delta$ are, up to natural isomorphism, mutually inverse functors from
$\mathbf{S}^{\mathfrak{P}}$ to the category sheaves over $\mathfrak{P}$. Thus, we
have:

\begin{theorem}
$\zip^\Delta$ and $\unzip^\Delta$ are a compression of $\mathbf{S}^{\mathfrak{P}}$.
\end{theorem}

In general, there is no canonical choice for admissible posets which automatically
glues over $\Delta$. However, below we will give a gluing procedure starting from
a family of admissible posets over $\Delta_{\max}$, which yields a set of admissible
posets together with a compression for any coherent $\Delta$-family.
For smooth toric varieties, the $\lcm$-lattices already will do the job:

\begin{proposition}
Assume that $X$ is a smooth toric variety, then the family $\zip^\sigma E^\sigma$,
with respect to the $\lcm$-lattices of the modules $E^\sigma$, glues over $\Delta$.
\end{proposition}

\begin{proof}
Without loss of generality assume that $\sigma$ has full dimension in $N_\R$. Let
$\tau < \sigma$ and for $m \in M$ denote $\bar{m}$ its class in $M / \tau^\bot_M$.
For any $m' \leq_\sigma m \in M$, an isomorphism $\chi^\sigma_{m', m} : E^\sigma_{m'}
\longrightarrow E^\sigma_m$ implies an isomorphism $\chi^\tau_{m', m} :
E^\tau_{m'} \longrightarrow E^\tau_m$, and for any $m \in M$, we have that $\bar{m}'
\in I^\tau(\bar{m})$ implies that there exists some $m'' \in \bar{m}$ with $m'' \in
I^\sigma(m)$. Therefore, if we denote $\p^\sigma$, $\p^\tau$ the $\lcm$-lattices of
$E^\sigma$ and $E^\tau$, respectively, we have a canonical contraction
$(\p^\sigma)_{\lessgtr_\sigma^\tau} \longrightarrow (\p^\tau)_{\lessgtr_\tau}$ given
by $\p^\sigma \ni m \mapsto A^\tau_E(\bar{m})$.
\end{proof}

\paragraph{Refining compressions.}
For some arbitrary choice of $\mathfrak{P}$, the category $\mathbf{S}^\mathfrak{P}$
in general contains not enough reflexive sheaves of rank one to construct global
resolutions. This is true even for the collection of
$\lcm$-lattices of \msh{E} over a smooth toric variety. The reason for this is that
relations of the module $E^\sigma$ which are encoded in the lattice $\p^\sigma$,
must no longer be present in the localization $E^\tau$ for $\tau < \sigma$, and thus
are ``contracted'' over $\p^\tau$. But if we construct a resolution over $U_\sigma$,
these relations still are present after we restrict to $U_\tau$, and thus are also
felt by neighbouring cones $\sigma'$ with $\tau \subset \sigma \cap \sigma'$. So,
in order to construct a resolution with respect to any admissible collection of
posets which glues over $\Delta$, we have to refine these posets in a way which
allows that any locally given reflexive sheaf $\sh{O}_{U_\sigma}(D)$ can be extended
to a suitable reflexive sheaf of rank one over $X$
 
\comment{
To construct a resolution of \msh{E} starting from the local data given by
compression with respect to admissible posets, one has to take into account that
every free resolution over some poset $\p^\sigma$ will interfer with the free
resolutions
over all the other $\p^\tau$'s; this is meant in the sense that for any $\n \in
\p^\sigma$ we have to extend the associated sheaf $\sh{O}_{U_\sigma}(\sum_{\rho}
n_\rho D_\rho)$, and correspondingly the map $\sh{O}_{U_\sigma} (\sum_{\rho} n_\rho
D_\rho) \longrightarrow \sh{E}\mid_{U_\sigma}$, to the whole of $X$. This corresponds
to extending in a suitable sense the map $F^\n \longrightarrow \zip^\sigma E^\sigma$
to maps over the other $\p^\sigma$. For this, the $\p^\sigma$ a priori are not fine
enough, as they do not take into account relations of the $E^\tau$ coming from
neighbouring $\p^\tau$, which are killed by restricting to $\p^{\sigma \cap \tau}$,
as these will be felt by globalized resolutions. So the first step to construct
global resolutions is to construct a refinements of the $\p^\sigma$, i.e. admissible
posets $\tilde{\p}^\sigma$ with $\p^\sigma \subset \tilde{\p}^\sigma$, such that
$\big(\tilde{\p}^{\sigma_1}\big)_{\lessgtr_{\sigma_1}^\tau} =
\big(\tilde{\p}^{\sigma_1}\big)_{\lessgtr_{\sigma_2}^\tau}$
for every $\tau, \sigma_1, \sigma_2 \in \Delta$, where $\tau = \sigma_1 \cap
\sigma_2$.
}

\begin{example}
\label{p1p1example}
Consider the toric surface $\mathbb{P}^1 \times \mathbb{P}^1$. The associated fan has
four rays $\rho^1, \dots, \rho^4$ and four maximal cones $\sigma^{12}, \sigma^{23},
\sigma^{34}, \sigma^{41}$, where $\sigma_{ij}$ is spanned by the rays $\rho^i$,
$\rho^j$. The associated semigroups $\sigma_M^{ij}$ are generated in $M \cong \Z^2$
by $\{(1, 0), (0, 1)\}$, $\{(0, 1), (-1, 0)\}$, $\{(-1, 0), (0, -1)\}$, $\{(0, -1),
(1, 0)\}$, respectively. We consider the sky\-scraper sheaf \msh{S} which has two
stalks at the orbits $\orb{\sigma^{12}}$ and $\orb{\sigma^{23}}$, respectively,
which are, as $k[\sigma_M^{ij}]$-modules, given by:
\begin{equation*}
\Gamma(U_{12}, \sh{S}) = k \cdot \chi\big((0, 0)\big), \quad
\Gamma(U_{23}, \sh{S}) = k \cdot \chi\big((-2, 2)\big).
\end{equation*}
Figure \ref{f-p1p1example} shows the four dual cones describing $\mathbb{P}^1 \times
\mathbb{P}^1$, slightly moved away from each other, and an indication of the
associated $\lcm$-lattices. The squares indicate the degrees $(0, 0)$ and $(-2, -2)$
where the stalks of \msh{S} sit, the light grey triangles indicate the anchor
elements of the two $\sigma$-families. 
\begin{figure}[htb]
\begin{center}
\includegraphics[width=8cm]{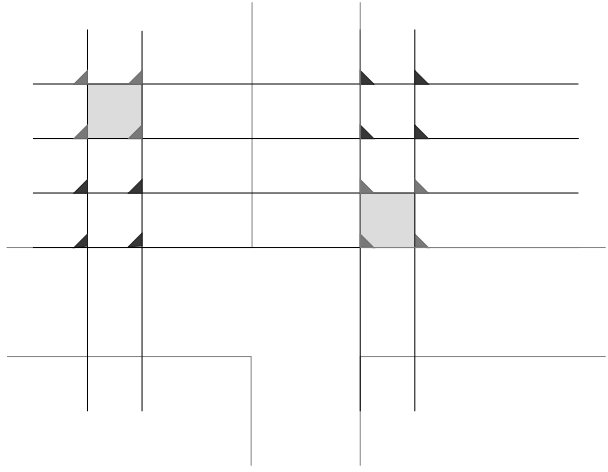}
\end{center}
\caption{Skyscraper sheaf over $\mathbb{P}^1 \times \mathbb{P}^1$.}\label{f-p1p1example}
\end{figure}
The dark grey triangles show the additional anchor elements which come from the
transition from one $\sigma$-family into another some of which have to enter a global
resolution. One possible resolution would be:
\begin{gather*}
0 \longrightarrow \sh{O}(-D_1 - D_2 - 2D_3) \oplus \sh{O}(-3D_2 - 3 D_3)
\longrightarrow \\
\sh{O}(-D_1 - 2 D_3) \oplus \sh{O}(-D_2-2D_3) \oplus \sh{O}(-2D_2 - 3 D_3) \oplus
\sh{O}(-3D_2 - 2 D_3) \longrightarrow \\ \sh{O}(-2D_3) \oplus
\sh{O}(-2 D_2 - 2 D_3) \longrightarrow
\sh{S} \longrightarrow 0
\end{gather*}
where we write $D_i$ instead of $D_{\rho_i}$. Note that the choice of other admissible
posets instead of the $\lcm$-lattices can lead to more convenient resolutions.
\end{example}

We consider any family of posets $\mathfrak{P} = \{\p^\sigma \mid \sigma
\in \Delta\}$, which glues over $\Delta$ and which is admissible with respect to
\msh{E}. We are going to construct a family of posets $\tilde{\mathfrak{P}} =
\{\tilde{\p}^\sigma \mid \sigma \in \Delta\}$ which glues over $\Delta$, is
admissible with respect to \msh{E}, and whose associated category
$\mathbf{S}^{\tilde{\mathfrak{P}}}$ has enough reflexive sheaves.
We start bottom-up and we consider $\p^\rho \subset \Z^{\rho(1)} \cong \Z$ for some
$\rho \in \rays$. In fact, $\p^\rho$ is a linear chain, i.e. a totally orderd
subset of $\Z$. For every $\sigma > \rho$, we consider
$\big(\p^\sigma\big)_{\lessgtr_\sigma^\rho}$ as a subset of $\Z^{\rho(1)}$, such that
the hooking $h_{\sigma\rho}$ becomes the natural inclusion $\p^\rho \subset 
\big(\p^\sigma\big)_{\lessgtr_\sigma^\rho}$ in $\Z^{\rho(1)}$. We set
\begin{equation*}
\tilde{\p}^\rho := \bigcup_{\rho < \sigma}\big(\p^\sigma\big)_{\lessgtr_\sigma^\rho}.
\end{equation*}
Now fix some $\sigma \in \Delta$ together with its admissible lattice
$\p^\sigma \subset \Z^{\sigma(1)}$. For every $\tau < \sigma$, we consider the
natural embedding $\Z^{\tau(1)} \hookrightarrow \Z^{\sigma(1)}$ which is induced
by the inclusion $\tau(1) \subset \sigma(1)$. In particular, every element $i \in
\tilde{\p}^\rho$ becomes an element of $\Z^{\sigma(1)}$ which is nonzero only at the
$\rho$th position. For every $m \in M$ there exists a unique anchor
element $A(m) \in \p^\sigma$. We refine now by setting:
\begin{equation*}
\tilde{A}(m) = \lcm \{i \in \tilde{\p}^\rho \mid i \leq \langle m,
n(\rho) \rangle\}_{\rho \in \sigma(1)}
\end{equation*}
and
\begin{equation*}
\tilde{\p}^\sigma := \{\tilde{A}(m) \mid m \in M\} \cup \p^\sigma,
\end{equation*}
where we observe that $A(m)$ is of the form
\begin{equation*}
A(m) = \big(\max\{i \in \tilde{\p}^\rho \mid i \leq A(m)_{\lessgtr_\sigma^\rho}\}
\mid \rho \in \sigma(1)\big)
\end{equation*}
and thus $\p^\sigma \subset \tilde{\p}^\sigma$
Clearly, $\tilde{\p}^\sigma$ is admissible with respect to $E^\sigma$, and we can
consider the compressions $\zip^{\tilde{P}^\sigma}$, $\unzip^{\tilde{P}^\sigma}$.
Using the identification of $(\p^\sigma)_{\lessgtr^\tau_\sigma}$ with its image
in $\Z^{\tau(1)}$, we have:

\begin{proposition}
For any $\tau \in \Delta$, $\tilde{\p}^\tau = \bigcup_{\tau < \sigma}
(\p^\sigma)_{\lessgtr^\tau_\sigma}$, where the union runs over all $\sigma
\in \Delta_{\max}$ with $\tau < \sigma$.
\end{proposition}

\begin{proof}
This follows because for any $\eta < \tau$, $\p^\eta = (\p^\eta)_{\lessgtr_\eta} \in
\Z^\tau$ is a subset of the image of $(\p^\sigma)_{\lessgtr_\sigma^\eta}$ in
$\Z^{\eta(1)}$.  Thus $\tilde{\p}^\rho = \bigcup_{\rho < \sigma}
(\p^\sigma)_{\lessgtr_\sigma^\rho}$ where the union runs over all maximal cones. Now
the proposition follows from $\p^\tau \subset (\p^\sigma)_{\lessgtr^\tau_\sigma}$
and by the generatedness of $\tilde{\p}^\sigma$ by $\p^\sigma$ and the
$\tilde{\p}^\rho$.
\end{proof}

By this proposition, we can conclude that the choice of any collection of admissible
posets leads to a collection of admissible posets which glue over $\Delta$:

\begin{corollary}
The family $\tilde{\p}^\tau$ is generated by the $\p^\sigma$, where $\tau$ runs over
$\Delta_{\max}$.
\end{corollary}

\begin{corollary}
$(\tilde{\p}^\sigma)_{\lessgtr_\sigma^\tau} = \tilde{\p}^\tau$ for all
$\tau < \sigma \in \Delta$.
\end{corollary}

By combining these two corollaries, we obtain:

\begin{proposition}
The family of sheaves $\zip^{\tilde{\p}^\sigma} E^\sigma$ glues over $\Delta$.
\end{proposition}

\paragraph{Global resolutions.}
Recall from section \ref{homext} that for every $\sigma \in \Delta$ we can construct
the extension module $EE^\sigma$ of $E^\sigma$ over the ring $S^\sigma$. In the
equivariant setting, the category of modules over $S^\sigma$ is equivalent to that
of the
ring $S_{x^{\hat{\sigma}}}$, and we can extend $EE^\sigma$ to a module over this
ring. By naturality of the construction, the $EE^\sigma$ glue to a sheaf $E\sh{E}$
over $\hat{X}$, and we obtain the $S$-module
$E\hat{E} := \Gamma(k^\rays, E\sh{E})$. We have the following properties for
$E\hat{E}$, which immediately follow from the corresponding properties of
proposition \ref{EEprop}:

\begin{proposition}
$E\hat{E}$ has the following properties:
\begin{enumerate}[(i)]
\item $E\hat{E}\breve{\ } \cong \sh{E}$.
\item If \msh{E} is coherent, then $E\hat{E}$ is finitely generated.
\item if \msh{E} is torsion free, then $E\hat{E}$ is torsion free.
\end{enumerate}
\end{proposition}

This way, a global resolution can be constructed as the descend of a
resolution of the $S$-module $E\hat{E}$ with respect to its $\lcm$-lattice. However,
there are more possibilites to resolve \msh{E} which use $E\hat{E}$ but do not
require the cost of computing the whole $\lcm$-lattice of $E\hat{E}$.
For this, we give a more precise picture of $E\hat{E}$.

For every $\sigma \in \Delta$ and every $\tau < \sigma$, denote $\pi_\sigma :
\Z^\rays \longrightarrow \Z^{\sigma(1)}$ and $\pi^\sigma_\tau: \Z^{\sigma(1)}
\longrightarrow \Z^{\tau(1)}$ the canonical projections. For any
$\n \in \Z^{\sigma(1)}$, $EE^\sigma_\n$ is defined to be the inverse limit
$E^\sigma(U_\n) := \underset{\leftarrow}{\lim} E^\sigma(U_m)$.
For any $\tau < \sigma \in \Delta$, there is the canonical map induced by
localization: $E^\sigma(U_\n) \longrightarrow E^\tau(U_{\pi^\sigma_\tau(\n)})$, and
for any $\n \in \weildivisors$, we obtain the directed system
\begin{equation*}
\xymatrix{
E\hat{E}_\n \ar[r]^{\pi_\sigma} \ar[rd]_{\pi_\tau} & EE^\sigma_{\pi_\sigma(\n)}
\ar[d]^{\pi^\sigma_\tau} \\
& EE^\tau_{\pi_\tau(\n)}
}
\end{equation*}
whose final object is the vector space $EE^0_{\pi_0(\n)}$.
The component $E\hat{E}_\n$ then has the universal property of the inverse limit of
this system:
\begin{equation*}
E\hat{E}_\n = \underset{\leftarrow}{\lim} E^\sigma_\n = \underset{\leftarrow}{\lim}
E^\sigma_m
\end{equation*}
where the latter limit runs over all $\sigma \in \Delta$ and the system of all $m
\in M$ such that $\n \leq_\sigma m$.

\begin{definition}
A {\em lift} $\tilde{\p}^\lambda$ of the collection $\tilde{\p}^\sigma$ is a
collection of injective,
order preserving maps $\lambda_\sigma: \tilde{\p}^\sigma \hookrightarrow \Z^\rays
\cup \hat{0}$ such that
\begin{enumerate}[(i)]
\item $\pi_\sigma\big(\lambda_\sigma(\n)\big) = \n$ for all $\n \in
\tilde {\p}^\sigma$;
\item $\lambda_\tau(\n) = \lcm\Big\{\lambda_\sigma \big((\pi^\tau_\sigma)^{-1}(\n)
\cap \tilde{\p}^\sigma\big) \mid \tau < \sigma\Big\}$
for every $\n \in \tilde{\p}^\tau$;
\item for all $\n \in \p^\sigma$ the composition $EE_{\lambda_\sigma(\n)}
\longrightarrow (EE_{\lambda_\sigma(\n)})_{\lessgtr^\sigma_\Delta} \longrightarrow
EE^\sigma_\n$ is surjective.
\end{enumerate}
We identify $\tilde{\p}^\lambda$ with the poset given by the image of the maps
$\lambda_\sigma$.
\end{definition}

There is, of course, no most natural choice for a lift $\lambda$, but a general
choice which always works, is:
\begin{equation*}
\lambda_\sigma(\n)_\rho =
\begin{cases}
n_\rho & \text{ if } \rho \in \sigma(1), \\
\max\{i \in \tilde{\p}^\rho\} & \text{ if } \tilde{\p}^\rho \neq \hat{0}, \\
0 & \text{ else}.
\end{cases}
\end{equation*}
In the case where \msh{E} is reflexive, it is possible to do a more efficient
general choice, as we will see in subsection \ref{reflext}. 
With respect to a lift $\lambda$, we can define the submodule $E\hat{E}_\lambda
\subset E\hat{E}$ as follows. For every $\sigma \in \Delta$ and all $\n \in
\tilde{\p}^\sigma$ we choose a subvector space $E'_{\lambda_\sigma(\n)} \subset
E\hat{E}_{\lambda_\sigma(\n)}$ such that the induced morphism
$E'_{\lambda_\sigma(\n)} \longrightarrow EE^\sigma_\n$ is surjective. Then we define
$E\hat{E}_\lambda$ to be the module generated by the $E'_{\lambda_\sigma(\n)}$.
Note that in spite we do not make explicit this choice in the notation, we always
assume it implicitly.

\begin{proposition}
$(E\hat{E}_\lambda)\breve{\ } \cong \sh{E}$.
\end{proposition}

\begin{proof}
As the homomorphism $EE_{\lambda_\sigma(\n)} \longrightarrow EE^\sigma_\n$ is
surjective, the map $(EE_{\lambda_\sigma(\n)})_{\lessgtr^\sigma_\Delta}
\longrightarrow EE^\sigma_\n$ becomes an isomorphism. Moreover, $\pi_\sigma \big(
\tilde{\p}^\sigma\big) = \tilde{\p}^\sigma$, so that the induced representation
on $\tilde{\p}^\sigma$ is the same as for $E\hat{E}$.
\end{proof}

So, we can also use $E\hat{E}_\lambda$ to construct resolutions:

\begin{corollary}
Every lift $\tilde{\p}^\lambda$ gives rise to a resolution of \msh{E}.
\end{corollary}

Instead of taking $E\hat{E}_\lambda$, we can also take directly the poset
$\tilde{\p}^\lambda$. Denote $i : \tilde{\p}^\lambda \hookrightarrow \Z^\rays$ the
canonical inclusion. We define:
\begin{equation*}
\zip^\lambda\sh{E} := i^* E\hat{E}_\lambda.
\end{equation*}
For any given representation $E$ of $\tilde{\p}^\lambda$, by the canonical projection
we obtain back the admissible poset $\tilde{\p}^\sigma$ as the image of
$\tilde{\p}^\lambda$ in $\Z^{\sigma(1)}$, together with the localization of the
representation $E$. These representations glue naturally over $\tilde{\p}^\sigma$,
and we can use them reconstruct the
module $E\hat{E}$. By taking the submodule generated in the degrees given by
$\tilde{\p}^\lambda$, we obtain $E\hat{E}_\lambda$. Using either $E\hat{E}$ or
$E\hat{E}_\lambda$, by sheafification we get back the sheaf \msh{E}. We denote this
procedure $\unzip^\lambda E$.

\begin{proposition}
The category of representations of $\tilde{\p}^\lambda$ is a full subcategory of
the category of sheaves which glue over the collection $\tilde{\p}^\sigma$.
\end{proposition}

\begin{proof}
We only remark that these categories in general can not be equivalent, as the lift
$\lambda^\sigma$ for every $\n \in \tilde{\p}^\sigma$ fixes the choice of free
representations of $\tilde{\p}^{\sigma'}$, $\sigma' \in \Delta$, which glue together
with the free representation of $\n$ over $\tilde{\p}^\sigma$.
\end{proof}

Using this correspondence, we obtain finally the finest class of global resolutions
for \msh{E}. For the representation $E^\lambda$ of $\tilde{\p}^\lambda$ for some
lift $\lambda$, we construct the free resolution in the category of
$\tilde{\p}^\lambda$-representations, and by unzipping we obtain:
\begin{equation*}
0 \longrightarrow \unzip^\lambda F_s \longrightarrow \cdots \longrightarrow
\unzip^\lambda F_0 \longrightarrow \sh{E} \longrightarrow 0.
\end{equation*}
However, nothing prevents us from taking a different lift $\lambda$ for every
syzygy of \msh{E}, and as we will see below, this will be quite natural for doing so
in the case of reflexive sheaves. So we can consider a sequence of lifts
$\lambda_0, \dots \lambda_s$ and a corresponding resolution:
\begin{equation*}
0 \longrightarrow  \unzip^{\lambda_s} F_s \longrightarrow \cdots \longrightarrow
\unzip^{\lambda_0} F_0 \longrightarrow \sh{E} \longrightarrow 0.
\end{equation*}
Here, the finiteness of the sequence follows that in every step we eliminate
minimal elements of the induced representations of the admissible posets
$\tilde{\p}^\sigma$, but the length $s$ finally may depend on the successive choice
of the lifts.

\begin{example}
\label{p1p1example2}
Consider the variety $\mathbb{P}^1 \times \mathbb{P}^1$ and skyscraper sheaf \msh{S}
similar to that of example \ref{p1p1example}, but this time with slightly different
gradings:
\begin{equation*}
\Gamma(U_{12}, \sh{S}) = k \cdot \chi(1, 1), \quad \Gamma(U_{23}, \sh{S}) = k \cdot
\chi(-1, 1).
\end{equation*}
Figure \ref{f-p1p1example2} shows the corresponding posets. Explicitly, we have:
\begin{align*}
\tilde{\p}^{12} & = \{\hat{0}, (1, 1), (2, 1), (1, 2), (2, 2)\} \\
\tilde{\p}^{23} & = \{\hat{0}, (1, 1), (2, 1), (1, 2), (2, 2)\} \\
\tilde{\p}^2 & = \{\hat{0}, 1, 2\}
\end{align*}
(for brevity, we suppress the ray $\rho_4$). We consider the lifts:
\begin{align*}
\lambda_{12}\big(\tilde{\p}^{12}\big) & = \{\hat{0}, (1, 1, 1), (2, 1, 2), (1, 2, 1),
(2, 2, 2)\} \\
\lambda_{12}\big(\tilde{\p}^{23}\big) & = \{\hat{0}, (1, 1, 1), (1, 2, 1), (2, 1, 2),
(2, 2, 2)\} \\
\lambda_{12}\big(\tilde{\p}^2\big) & = \{\hat{0}, (2, 1, 2), (2, 2, 2)\}.
\end{align*}
\begin{figure}[htb]
\begin{center}
\includegraphics[width=8cm]{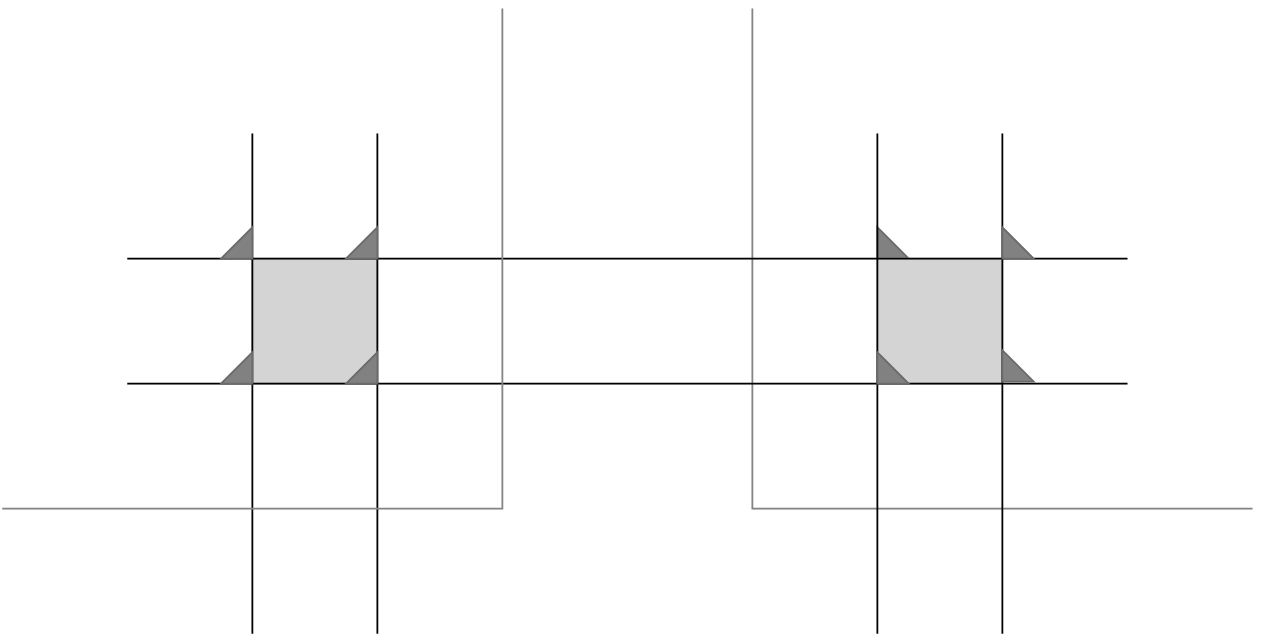}
\end{center}
\caption{Another skyscraper sheaf over $\mathbb{P}^1 \times \mathbb{P}^1$.}\label{f-p1p1example2}
\end{figure}
The sheaf $E\hat{E}$ is given by
\begin{equation*}
E\hat{E}_\n \cong
\begin{cases}
k^2 & \text{ if } \n = (1, 1, 1) \\
0 & \text{ else}.
\end{cases}
\end{equation*}
Choosing the one-dimensional diagonal $k \subset E \hat{E}_{(1,1,1)}$, we obtain as
resolution:
\begin{gather*}
0 \longrightarrow \sh{O}(- 2 D_1 - 2 D_2 - 2 D_3) \longrightarrow
\sh{O}(- 2 D_1 - D_2 - 2 D_3) \oplus \sh{O}(- D_1 - 2 D_2 - D_3) \\ \longrightarrow
\sh{O}(-D_1 - D_2 - D_3)
\longrightarrow \sh{S} \longrightarrow 0
\end{gather*}
\end{example}

\section{Reflexive Sheaves and Vector Space Arrangements}
\label{reflexivesheaves}

\subsection{Reflexive sheaves and their canonical admissible posets}

For an equivariant reflexive sheaf over a toric variety, i.e. a sheaf \msh{E} which is
isomorphic to its bidual, $\sh{E} \cong \sh{E}\check{\ }\check{\ }$, the associated
$\Delta$-family has a quite efficient representation. To every equivariant
coherent sheaf \msh{E} over $U_\sigma$, one can associate a limit vector space
$\mathbf{E}^\sigma := \underset{\rightarrow}{\lim} E_m^\sigma$, and by the gluing
of the $E^\sigma$ over the collection of posets $(M, \leq_\sigma)$, there is a
functorial isomorphism $\mathbf{E}^\sigma \rightarrow \mathbf{E}^0 =: \mathbf{E}$,
where $0$ denotes the zero cone in $\Delta$, and moreover, $\dim \mathbf{E} = \rk
\sh{E}$. As explained in detail in
\cite{perling1}, section 5 (see also \cite{Kly90}, \cite{Kly91}), every
equivariant reflexive sheaf \msh{E} is determined by a set of filtrations
\begin{equation*}
\cdots \subset E^\rho(i) \subset E^\rho(i + 1) \subset \cdots \subset \mathbf{E}
\end{equation*}
for every ray $\rho \in \rays$. These filtrations must be {\em full}, i.e.
$E^\rho(i) = 0$ for very small $i$, and $E^\rho(i) = \mathbf{E}$
for $i$ very large. The corresponding $\sigma$-families then can be constructed from
these filtrations by setting
\begin{equation*}
E^\sigma_m = \bigcap_{\rho \in \sigma(1)} E^\rho\big(\langle m, n(\rho) \rangle\big).
\end{equation*}
In fact, this construction establishes an equivalence of categories between
equivariant reflexive sheaves and vector spaces with full filtrations.
The morphisms in the latter category are vector space homomorphisms which are
compatible with the filtrations in the $\Delta$-family sense (\cite{perling1},
Theorem 5.29).

Consider a reflexive module $E^\sigma$ over the ring $k[\sigma_M]$, where without
loss of generality we assume that $\sigma$ has full dimension in $N_\R$. To any
such module there is associated the {\em subvector space arrangement} $\{E^\sigma_m
\mid m \in M\}$ in the limit vector space $\mathbf{E}$, where $E_m^\sigma =
\bigcup_{\rho \in \sigma(1)} E^\rho\big(\langle m, n(\rho) \rangle\big)$. This
arrangement in a natural way is a poset, where the partial order is given by
inclusion. We will show that we can embed this poset into $\Z^{\sigma(1)}$ such
that it becomes an admissible poset for $E^\sigma$.

\begin{definition}
For every $m \in M$, we define $\kappa_\rho(m) = \min \{i \in \Z \mid E_m^\sigma
\subset E^\rho(i)\}$ and the {\em anchor} of $m$ by:
\begin{equation*}
A(m) = \big(\kappa^\rho_m \mid \rho \in \sigma(1)\big) \in \Z^{\sigma(1)}.
\end{equation*}
We denote $\p_{E^\sigma}$ the subposet $\{A(m) \mid m \in M\}$ of
$\Z^{\sigma(1)}$.
\end{definition}

\begin{proposition}
\label{admissibleproof}
$\p_{E^\sigma}$ is admissible with respect to $E^\sigma$.
\end{proposition}

\begin{proof}
First, clearly, $A(m) \leq m$ for all $m \in M$. Now assume that $A(m') \leq m$ for
some $m' \in M$. $A(m') \leq m$ implies that $E^\sigma_{m'} \subset E^\sigma_m$,
and thus $A(m') \leq A(m)$. Now, by definition $\bigcap_{\rho \in \sigma(1)} E^\rho
\big(n_\rho) = E^\sigma_m$ for all $m \in T_\n$ for some $\n \in \p_{E^\sigma}$.
\end{proof}

\begin{definition}
We call $\p_{E^\sigma}$ the {\em canonical admissible poset} of $E^\sigma$.
\end{definition}

An important fact for understanding the structure of reflexive modules is the
following

\begin{lemma}
\label{reflkeylemma}
Let $\p_{E^\sigma}$ be the canonical admissible poset of $E^\sigma$. Then $E^\sigma_m
\subset E^\sigma_{m'}$ iff $A(m) \leq_\sigma A(m')$. Moreover, $E^\sigma_m =
E^\sigma_{m'}$  iff $A(m) = A(m')$.
\end{lemma}

\begin{proof}
Assume first that $E^\sigma_m \subset E^\sigma_{m'}$. Then for every $\rho \in
\sigma(1)$ it follows that $\min\{i \mid E^\sigma_m \subset E^\rho(i)\} \leq \min\{i
\mid E^\sigma_{m'} \subset E^\rho(i)\}$, and thus $A(m) \leq_\sigma A(m')$. In the
other direction, denote $\n := A(m)$, $\n' := A(m')$, then $n_\rho \leq n'_\rho$ for
every $\rho \in \sigma(1)$ and $E^\rho(n_\rho) \subseteq E^\rho(n'_\rho)$, and thus
$E^\sigma_m \subset E^\sigma_{m'}$.
\end{proof}

\begin{proposition}
If $U_\sigma$ is smooth, then, as a poset, the vector space arrangement associated
to $E^\sigma$ is isomorphic to its $\lcm$-lattice.
\end{proposition}

\begin{proof}
Because $U_\sigma$ is smooth, for every $m \in M$, the anchor element $A(m)$ is an
element of $M$, and we conclude from the proof of proposition \ref{admissibleproof},
that $A(m)$ is the unique member of $I(m)$. For any two $A(m) \neq A(m')$, the vector
spaces $E^\sigma_m$ and $E^\sigma_{m'}$ do not coincide, and thus the vector space
$E^\sigma_{m''}$, where $m'' = \lcm \{m, m'\}$, contains at least the sum $E^\sigma_m
+ E^\sigma_{m'}$. Moreover, we have that $E^\sigma_{m''} = \bigcap_{\rho \in
\sigma(1)} E^\rho\big(\langle m'', n(\rho) \rangle\big)$, where for every $\rho \in
\sigma(1)$ $E^\rho\big(\langle m'', n(\rho) \rangle\big)$ contains $E^\sigma_m$ and
$E^\sigma_{m'}$, and thus $\langle m'', n(\rho) \rangle \geq \max\{\langle m, n(\rho)
\rangle, \langle m', n(\rho) \rangle\}$. So $m''$ is the minimal element of $M$
with respect to the partial order $\sigma_M$, such that $E^\sigma_{m''}$ contains
both, $E^\sigma_m$ and $E^\sigma_{m'}$.
\end{proof}

\begin{example}
\label{intersectionimprove}
We give an example which shows that the choice of another admissible poset instead of
the canonical one can improve the resolution.
Consider the subsemigroup $\sigma_M$ of $\Z^2$ which is generated by $(1, 0)$,
$(1, 1)$ and $(1, 2)$; the corresponding cone $\sigma$ has two rays $\rho_1, \rho_2$
with primitive elements $n(\rho_1) = (2, 1)$, $n(\rho_2) = (0, 1)$. Let $\mathbf{E}
\cong k^3$ and consider the filtrations
\begin{equation*}
E^{\rho_1}(i) =
\begin{cases}
0 & \text{ for } i < 0 \\
E_1 & \text{ for } i = 0 \\
\mathbf{E} & \text{ for } i > 0
\end{cases}\qquad
E^{\rho_2}(i) =
\begin{cases}
0 & \text{ for } i < 1 \\
E_2 & \text{ for } i = 1 \\
\mathbf{E} & \text{ for } i > 1.
\end{cases}
\end{equation*}
with $\dim E_i = 2$ and the $E_i$ in general position.
The corresponding canonical admissible poset is
$\p = \{\hat{0}, (0, 2), (1, 1), (1, 2)\}$ and it leads to the resolution
\begin{equation*}
0 \longrightarrow S_{(2, 2)} \longrightarrow S_{(1, 1)}^2 \oplus S_{(0, 2)}^2
\longrightarrow E \longrightarrow 0
\end{equation*}
\begin{figure}[htb]
\begin{center}
\includegraphics[width=8cm]{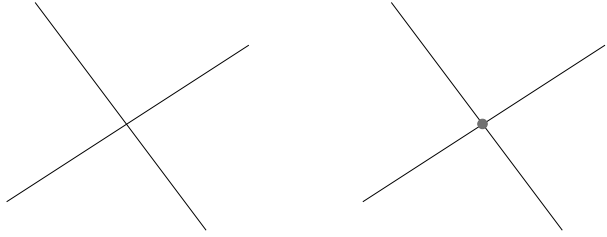}
\end{center}
\caption{Canonical admissible poset and the poset generated by its intersections}\label{f-intersectionimprove}
\end{figure}
If we choose instead the poset $\p'= \{\hat{0}, (0, 1), (0, 2), (1, 1), (1, 2)\}$,
the associated representation of $\p'$ maps $(0, 1)$ to the the subvector space
$E_1 \cap E_2$ of $\mathbf{E}$. The corresponding vector space arrangements are shown
as linear configurations in $\mathbb{P}\mathbf{E} \cong \mathbb{P}^2$ in figure
\ref{f-intersectionimprove}. The grey dot in the right figure denotes the intersection
$E_1 \cap E_2$. The corresponding resolution becomes:
\begin{equation*}
0 \longrightarrow S_{(0, 1)} \oplus S_{(0, 2)} \oplus S_{(1, 1)} \longrightarrow E
\longrightarrow 0,
\end{equation*}
i.e. $E$ splits into a direct sum of reflexive sheaves of rank one.
\end{example}

\subsection{Extensions to the homogeneous coordinate ring}
\label{reflext}

We first investigate the structure of the module $E\hat{E}$ where \msh{E} is
reflexive.
For this, we first consider the module $EE^\sigma$ for any $\sigma \in \Delta$. Its
determination is a straightforward computation:

\begin{proposition}
Let $E^\sigma$ be a reflexive $k[\sigma_M]$-module given by filtrations $E^\rho(i)$.
Then its extension is given by:
\begin{equation*}
EE^\sigma_\n = \bigcap_{\rho \in \sigma(1)} E^\rho(n_\rho).
\end{equation*}
\end{proposition}

\begin{proof}
We have $EE^\sigma_\n = \underset{\leftarrow}{\lim} E_m^\sigma$, where the limit runs
over all $\n \leq m$. As all morphisms $\chi^\sigma_{m, m'}$ are injective, this
direct limit immediately translates into an intersection in $\mathbf{E}^\sigma$:
\begin{align*}
\underset{\leftarrow}{\lim} E_m^\sigma & = \bigcap_{\n \leq m} E_m^\sigma \\
& = \bigcap_{\n \leq m} \bigcap_{\rho \in \sigma(1)} E^\rho\big(\langle m, n(\rho)
\rangle \big).
\end{align*}
It is always possible to find $m \in M$ for some $\tau \in \sigma(1)$ such that
$\langle m, n(\tau) \rangle = n_\tau$ and $\langle m, n(\rho) \rangle >> 0$ for any
$\tau \neq \rho$, such that $\bigcap_{\rho \in \sigma(1)} E^\rho\big(\langle m,
n(\rho) \rangle\big) = E^\tau(n_\tau)$. Thus we obtain $\bigcap_{\rho \in \sigma(1)}
E^\rho(n_\rho) \subset EE^\sigma_\n \subset \bigcap_{\rho \in \sigma(1)}
E^\rho(n_\rho)$ and the proposition follows.
\end{proof}

So the module $E\hat{E}^\sigma$ can explicitly be described by the filtrations for
\msh{E} and in fact, it is a reflexive module. To describe its filtrations more
explicitly, we use the quotient representation $\pi : k^{\sigma(1)} \longrightarrow
U_\sigma$.
For each $\rho \in \sigma(1)$, the restriction of \msh{E} to $U_\rho$ is a locally
free sheaf and thus if we restrict $\pi$ to $U_{\hat{\rho}}$, the pullback
\begin{equation*}
\hat{\sh{E}}^{\hat{\rho}} :=
(\pi\vert_{U_{\hat{\rho}}})^*\sh{E}\vert_{U_{\rho}}
\end{equation*}
is locally free over $U_{\hat{\rho}}$. To determine the filtration associated to
$\hat{\sh{E}}^{\hat{\rho}} $, consider the injective map
\begin{equation*}
\alpha_\rho: M / \rho^\bot_M \longrightarrow \Z^{\rho(1)}.
\end{equation*}  
Then every element $i \in \weildivisors / \hat{\rho}^\bot_{\hat{M}}$ lies in a unique
intervall $\alpha_\rho(j) \leq i < \alpha_\rho(j + 1)$ for some $j \in  M /
\rho^\bot_M \cong \mathbb{Z}$. $\hat{\sh{E}}^{\hat{\rho}}$ then can be described
by a filtration of $\mathbf{E}$, which is given by
\begin{equation*}
E\hat{E}^{\hat{\rho}}(i) = E^\rho(j) \text{ for } \alpha_\rho(j) \leq i
< \alpha_\rho(j + 1).
\end{equation*}
The reflexive $S$-module defined by set of filtrations $E\hat{E}^{\hat{\rho}}(i)$
for every $\rho \in \rays$ then can be identified with $E\hat{E}$.

\begin{proposition}
Let \msh{E} be a reflexive sheaf, then there is an isomorphism
$\hat{E}\check{\ }\check{\ } \cong E\hat{E}$.
\end{proposition}

\subsection{Resolutions for vector space arrangements and reflexive equivariant
sheaves}

\paragraph{The affine case.}
First we consider resolutions for a reflexive $M$-graded module $E^\sigma$ over
$k[\sigma_M]$ with filtrations $E^\rho(i)$ for $\rho \in \sigma(1)$.  Revisiting the
resolution process of proposition \ref{repres} for the corresponding representation
of the canonical admissible poset $\p_{E^\sigma}$, we find by lemma
\ref{reflkeylemma} that for any $\n \in \p_{E^\sigma}$, the vector space
$E^\sigma_{< \n}$ is the subvector space of $E^\sigma_\n$ which is spanned by all its
{\em sub}vector spaces in the arrangement $\p_{E^\sigma}$. We have the first step of
its resolution
\begin{equation*}
0 \longrightarrow K_0 \longrightarrow F_0 \longrightarrow E^\sigma \longrightarrow 0
\end{equation*}
such that $F_0$ is a reflexive module $F_0 \cong \bigoplus_{\n \in \p^\sigma}
S_{(\n)}^{f_\n}$
which is defined by filtrations $F^\rho(i)$ in a limit vector space $\mathbf{F}$,
defining a vector space arrangement $\mathcal{Q} := \{F_m \mid m \in M\}$.

\begin{proposition}
The poset underlying the vector space arrangement $\mathcal{Q}$ is isomorphic to
$\p_{E^\sigma}$.
\end{proposition}

\begin{proof}
The dimension of the vector space $F_{0, \n}$ is given by the number of $\n'
\leq_\sigma \n$; by lemma \ref{reflkeylemma} we have that $E^\sigma_m \subsetneq
E^\sigma_{m'}$ iff $A(m) < A(m')$, and thus the number of $\n'' \in \p$ for which
$E^\sigma_m$ has positive free dimension and which $\n'' \leq A(m)$ is smaller than
the number of such elements with $\n'' \leq A(m')$.
\end{proof}

The kernel $K_0$ is a reflexive module, given by filtrations $K^\rho(i) =
\operatorname{ker}(F^\rho(i)\rightarrow E\rho(i))$ of the kernel vector space
$\mathbf{K} = \operatorname{ker}(\mathbf{F} \rightarrow \mathbf{E})$. However,
the canonical admissible poset of $K_0$ is no longer isomorphic to $\p_{E^\sigma}$,
but we have the following:

\begin{proposition}
The canonical admissible poset of $K_0$ is a contraction of $\p_{E^\sigma}$.
\end{proposition}

\begin{proof}
We define the retraction morphism $r : \p_{E^\sigma} \longrightarrow \p_{K_0}$ by
mapping $A_{E^\sigma}(m)$ to $A_{K_0}(m)$ for all $m \in M$. For any $E^\sigma_m
\subset E^\sigma_{m'}$ we have $K_{0, m} \subset K_{0, m'}$, and thus
$r\big(U(A_{E^\sigma}(m))\big) \subset U(A_{K_0}(m))$. The other inclusion follows
because $\p_{K_0}$ is admissible for $K_0$. On the other hand, let $\n \in
\p_{K_0}$, then $\n \leq \n'$ for every $\n' \in r^{-1}\big(U(\n)\big)$, and
$\n \in  r^{-1}\big(U(\n)\big)$, thus $r^{-1}\big(U(\n)\big) = U(\n)$ in
$\p_{E^\sigma}$.
\end{proof}

By \ref{contractionliftres} this in particular implies that we can iterate and
the resolution of the vector space arrangement $\p_{E^\sigma}$ is equivalent to
a resolution of $E^\sigma$. We have:
\begin{equation*}
0 \longrightarrow F_s \longrightarrow \cdots \longrightarrow F_0 \longrightarrow
E^\sigma \longrightarrow 0
\end{equation*}
where $F_i \cong \bigoplus_{\n \in \p_{E^\sigma}} S_{(\n)}^{f^i_\n}$.

The shape of the resolution can be changed by chosing another admissible poset for
$E^\sigma$. This in turn is equivalent to adding {\em any} set of intersections of
vector spaces in $\p_{E^\sigma}$. To see this, we pass to the module $EE^\sigma$.
The arrangement of this module is complete with respect to intersections, and
every anchor element of the canonical admissible poset of $E^\sigma$ is by definition
an anchor element of the $\lcm$-lattice of $EE^\sigma$. In particular, for every
$\n \in \cL_{EE^\sigma}$ with $\n \leq m$, we have $\n \leq A_{E^\sigma}(\n)$
by lemma \ref{reflkeylemma}, so that condition (\ref{admissibledefi}) of definition
\ref{admissibledef} is fulfilled. Moreover, as $T_\n$ is empty if $\n$ is not from
$\p_{E^\sigma}$, condition (\ref{admissibledefii}) is trivially fulfilled.

\paragraph{The global case.}
Now we assume that \msh{E} is a reflexive sheaf over an arbitrary toric variety $X$,
represented by filtrations $E^\rho(i)$ of some vector space $\mathbf{E}$ for every
$\rho \in \rays$. We denote
$\p^\sigma$ the canonical admissible posets for every $E^\sigma$. To make contact
with the formalism of section \ref{deltaglobres}, we first consider the refinements
$\tilde{\p}^\sigma$.

\begin{lemma}
$\p^\sigma$ is a contraction of $\tilde{\p}^\sigma$ for every $\sigma \in \Delta$.
\end{lemma}

\begin{proof}
For every $\rho \in \rays$, the canonical admissible poset $\p^\rho$ is given by
$\hat{0}$ and some sequence $i^\rho_1 < \dots < i^\rho_{k_\rho}$ in $\Z$, where
$k_\rho < \rk \sh{E}$, such that $E^\rho(i) = E^\rho(i + j)$ for $j \geq 0$ if and
only if there exists no $i_{p}^\rho$ for some $p \in \{1, \dots, k_\rho\}$
such that $i < i^\rho_p \leq i + j$.
For every $\rho \in \rays$ and every $\rho < \sigma$, we have
$(\p^\sigma)_{\lessgtr_\sigma^\rho} = \p^\rho$, and thus $\p^\rho = \tilde{\p}^\rho$.
Recall that $\tilde{A}^\sigma$ was defined as the least common multiple of the
elements $\max\{i \in \tilde{\p}^\rho \mid i \leq \langle m, n(\rho) \rangle\}$,
where $\tilde{\p}$ is considered as subset of $\Z^{\sigma(1)}$ via the canonical
embedding $\Z^\rho \hookrightarrow \Z^{\sigma(1)}$. Denote $r : \tilde{\p}^\sigma
\longrightarrow \p^\sigma$, mapping the anchor $\tilde{A}^\sigma(m)$ to
$A^\sigma(m)$. Clearly, $r$ is surjective. Then for any $\tilde{A}^\sigma(m) \in
\tilde{\p}^\sigma$, the image of $U\big(\tilde{A}^\sigma(m)\big)$ is
$U\big(A^\sigma(m)\big)$. For any $\n \in \p^\sigma$, $r^{-1}(\n) = \n$, so
$r^{-1}\big(U(\n)\big) = U(\n)$ (the latter as an open subset of $\tilde{\p}^\sigma$,
and the lemma follows.
\end{proof}

For resolving \msh{E}, we now must define a lift $\lambda$ of the collection
$\tilde{\p}^\sigma$ to \mweildivisors. For every $\sigma \in \Delta$, we define
$\lambda_\sigma : \tilde{\p}^\sigma \longrightarrow \weildivisors$ by
\begin{equation*}
\big(\lambda_\sigma(\n)\big)_\rho =
\begin{cases}
\min\{i \mid E^\sigma_\n \subset E^\rho(i)\} & \text{ for } \rho \in \rays \setminus
\sigma(1) \\
n_\rho & \text{ for } \rho \in \sigma(1).
\end{cases}
\end{equation*}

\begin{proposition}
The collection $\lambda_\sigma$ is a lift of $\tilde{\p}^\sigma$.
\end{proposition}

\begin{proof}
By definition, $(\pi_\sigma \circ \lambda_\sigma)(\n) = \n$ for every $\n \in
\tilde{\p}^\sigma$. We show that $\lambda_\tau(\n) = \lcm\Big\{\lambda_\sigma
\big((\pi^\tau_\sigma)^{-1}(\n) \cap \tilde{\p}^\sigma\big) \mid \tau < \sigma\Big\}$
for every $\tilde{\p}^\tau$.
For this, observe that $E\hat{E}_{\lambda_\sigma(\n)} = E^\sigma_\n$, because
\begin{align*}
E^\sigma_\n & = \bigcap_{\rho \in \sigma(1)} E^\rho(n_\rho) 
\subset E\hat{E}_{\lambda_\sigma(\n)} 
= \bigcap_{\rho \in \rays} E^\rho( \lambda_\sigma(\n)_\rho) \\
&  = E^\sigma_\n \cap \big(\bigcap_{\rho \in \rays \setminus \sigma(1)}
E^\rho(\lambda_\sigma(\n)_\rho)\big) \subset E^\sigma_\n.
\end{align*}
\end{proof}

Now, the lift $\lambda$ gives rise to a subarrangement of the subvector space
arrangement of the arrangement associated to $E\hat{E}$, which is given by the
union of arrangements in $\mathbf{E}$:
\begin{equation*}
\p^\Delta := \bigcup_{\sigma \in \Delta} \p^\sigma = \Big\{\bigcap_{\rho \in
\sigma(1)} E^\rho\big(\langle m, n(\rho) \rangle\big) \mid \sigma \in \Delta, m \in
M \Big\} = \bigcup_{\sigma \in \Delta} \{\lambda_\sigma(\n) \mid \n \in \p^\sigma\}
\end{equation*}
The first step $0 \rightarrow \sh{K}_0 \rightarrow \sh{F}_0 \rightarrow \sh{E}
\rightarrow 0$ of the global resolution of $\mathbf{E}$ then is given by the sheaf
\begin{equation*}
\sh{F}_0 \cong \bigoplus_{\n \in \p^\Delta} \sh{O}\big(D_{\lambda(\n)}\big)^{f^0_\n},
\end{equation*}
where $f^0_\n$ is the free dimension of the vector space $E_\n$. By iteration, we
get a free resolution, which at the same time is a resolution of the vector
space arrangement $\p^\Delta$.
Note that this resolution coincides with the resolution of the module
$E\hat{E}_\lambda$ over $S$. The global resolution of \msh{E}
constructed using $E\hat{E}$ is given by the
minimal resolution given by the vector space arrangement in $\mathbf{E}$ which is
generated by {\em all} intersections of the vector spaces $E^\rho(i)$.

\subsection{Resolutions of Cohen-Macaulay modules}
\label{cmmodules}

Let $E$ be a (maximal) Cohen-Macaulay module over $k[\sigma_M]$, where $\sigma$
has full dimension in $N_\R$. We show that our resolutions behave well in the sense
that the maximal length of regular sequences does not decrease. We follow
\cite{brunsherzog} \S 1.5, and say that the graded module $E$ is Cohen-Macaulay if
$\operatorname{grade}_\mathfrak{m} E = \dim k[\sigma_M]$, where $\mathfrak{m}$ is the
maximal homogeneous ideal of $k[\sigma_M]$ which is generated by all non-unit
monomials.

\begin{theorem}
\label{CMresolution}
Let $E$ be an $M$-graded Cohen-Macaulay module over $k[\sigma_M]$ and consider the
resolution
\begin{equation*}
0 \longrightarrow F_s \longrightarrow \cdots \longrightarrow F_0 \longrightarrow E
\longrightarrow 0
\end{equation*}
corresponding to the canonical admissible poset of $E$. Then every $F_i$ is a
direct sum of Cohen-Macaulay modules of rank one.
\end{theorem}

\begin{proof}
We need only to consider the first step of the resolution $0 \rightarrow K_0
\rightarrow F_0 \rightarrow E \rightarrow 0$, as $K_0$ will be Cohen-Macaulay if
$F_0$ and $E$ are Cohen-Macaulay; the result then follows by induction. If we
restrict the surjection from $F_0$ to $E$ to a direct summand of rank one $R$ of
$F_0$, we necessarily obtain an injection $0 \rightarrow R \rightarrow E$.
We show that any $E$-regular sequence by construction also is a $R$-regular sequence.
Let $x_1, \dots, x_r$ be a $E$-regular sequence and denote $\mathbf{x}_i$ the ideal
generated by $x_1, \dots, x_i$, for $1 \leq i \leq r$. We consider the diagram
\begin{equation*}
\xymatrix{
& 0 \ar[d] & 0 \ar[d] & & \\
0 \ar[r] & \mathbf{x}_i R \ar[r] \ar[d] & \mathbf{x}_i E \ar[r] \ar[d] & \mathbf{x}_i
E / \mathbf{x}_i R \ar[r] \ar[d]^\alpha & 0 \\
0 \ar[r] & R \ar[r] \ar[d] & E \ar[r] \ar[d] & E / R \ar[r] & 0 \\
& R / \mathbf{x}_i R \ar[r]^\beta \ar[d] & E / \mathbf{x}_i E \ar[d] & & \\
& 0 & 0 & &
}
\end{equation*}
If $\alpha$ is injective, then also $\beta$ is injective, and the element $x_{i + 1}$
is a nonzero divisor of $R / \mathbf{x}_i R$, as it is a nonzero divisor of $E /
\mathbf{x}_i E$. To show that $\alpha$ is injective, we show that there exists no
$e_1, \dots, e_i \in E$ such that $y := \sum_{j = 1}^i x_j e_j$ is in $R$ but not in
$\mathbf{x}_i R$. This sum decomposes into homogeneous summands $y = \sum_{m \in M}
y_m$ where $y_m = \sum_{j = 1}^i \sum_{m' \in M} x_{j, m'} \cdot e_{j, m - m'}$.
If we write $x_{j, m'} = a_{j, m'} \chi(m')$, this sum can be written as
$\sum_{j = 1}^i \sum_{m' \in M} a_{j, m'} \chi(m') \cdot e_{j, m - m'}$. Now we
split the set $\{m' \in M \mid e_{j, m - m'} \neq 0\} = U_j \coprod V_j$, where
$U_j = \{m' \mid R_{m - m'} \neq 0\}$. By construction of the inclusion of $R$ in
$E$, there does not exist any $m'' \in V_j$ such that $E_{m''}$ contains a one
dimensional subvector space whose image in $\mathbf{E}$ coincides with the image of
$R$.  For any $m' \in V_j$, the elements
$\chi(m') \cdot e_{j, m - m'}$ must be contained in the subvector space $F_m$ spanned
by all $E_{m''}$ with $m'' < m$, and writing the equations modulo $F_m$, we can
replace every $e_j$ by some $f_j$ such that $f_{j, m - m'} = 0$ if $m' \in V_j$
and $\sum_j x_j f_j = y$. Thus we have for every $m$ the equation $x_m =
\sum_{j = 1}^i \sum_{m' \in U_j} a_{j, m'} \chi(m') \cdot f_{j, m - m'}$.
For $m' \in U_j$, we can project every
$f_{j, m - m'}$ to some appropriate $r_{j, m - m'} \in R_{m - m'}$, such that
$x_m = \sum_{j = 1}^i \sum_{m' \in U_j} a_{j, m'} \chi(m') \cdot r_{j, m - m'}
$.
Therefore, we have $x_m \in \mathbf{x}_i R$, from which follows that $\alpha$ is
injective.
\comment{
By construction
of the inclusion of $R$ in $E$, there does not exist any $m'' \leq_\sigma m$ such
that $E_{m''}$ contains a one dimensional subvector space whose image in $\mathbf{E}$
coincides with the image of $R$. Thus the summands in the sum over $V_j$ are all
contained in a proper subvector space of $E_m$ which does not contain $R_m$. Hence
we have a sum of two vectors in $E_m$ which lies in $R_m$, where one of the summands
lies in $R_m$, and the other outside, and so the second summand must vanish.
Therefore, we have $e_{j, m - m'} = r_{j, m - m'}$ for every nonzero $e_{j, m - m'}$,
and $e_1, \dots, e_i \in R$, from which follows that $\alpha$ is injective.}
\comment{
$x = \sum_{j = 1}^t \chi(m_j) f_j$, where the $f_j$ are homogeneous elements of $E$.
The image of $R$ in $\mathbf{E}$ spans a one-dimensional subvector space $\langle
R \rangle$ of $\mathbf{E}$, and for $x$ to be in $R$, it is a necessary condition
that the image of every summand $\chi(m_j) f_j$ is contained in $\langle R \rangle$.
The $\chi(m_j)$ act as identity homomorphism on the vector space $\mathbf{E}$ and
thus on its subvector spaces. So, it follows that every $f_j$ must be in $\langle R
\rangle$. By construction, the submodule $R$ of $E$ is the reflexive module
associated to an anchor element $\n \in \Z^{\sigma(1)}$, and by this it lies outside
of the span of all subvector spaces of $\bigcap_{\rho \in \sigma(1)} E^\rho(n_\rho)$,
hence the $f_j$ must be contained in $R$, and thus $x \in \mathbf{x}_i R$,
}
\end{proof}

\begin{corollary}[from proof of theorem \ref{CMresolution}]
Let $E$ be any reflexive $k[\sigma_M]$-module and $F_i$ as in the theorem, then
$\operatorname{grade}_\mathfrak{m} F_i \geq \operatorname{grade}_\mathfrak{m} E$ for
all $0 \leq i \leq s$.
\end{corollary}

\subsection{Reflexive models for vector space arrangements}
\label{reflexivemodels}

In this subsection we want to make a few remarks on how resolutions of vector space
arrangements can efficiently be constructed by passing to appropriate reflexive
modules over the polynomial ring. The point here is resolutions of such modules are
a standard task for many computer algebra systems. However, to make use of such
systems, one has to construct appropriate input data from the arrangement.
Let $\mathcal{V}$ be a subvector space arrangement of some vector space $\mathbf{V}$.
We make two assumptions on $\mathcal{V}$; the first is that $\mathcal{V}$ is complete
with respect to intersections, that is, for any subset $W_1, \dots, W_r \in
\mathcal{V}$, the intersection $W_1 \cap \dots \cap W_n$ is also in $\mathcal{V}$.
The second assumption is that the input data for $\mathcal{V}$ is given by a set of
vectors $v^W_1, \dots, v^W_{i_W}$ such that $W$ is the span over $k$ of all $v^V_i$
where $V \subset W$ and $i = 1, \dots, i_V$. Moreover, we assume that this set is
irredundant, i.e. $i_W = \codim_W \sum_{V \subsetneq W} V$.
Using this input data, the first step
\begin{equation*}
0 \longrightarrow K_0 \longrightarrow F_0 \overset{M}{\longrightarrow} \mathbf{V}
\longrightarrow 0
\end{equation*}
of the resolution of $\mathcal{V}$ is nearly tautological. Assume that we have chosen
a basis for $\mathbf{V}$, then $F_0$ is given by a basis $e^W_i$, $i = 1, \dots, i_W$
in one-to-one correspondence to the vectors $v^W_i$, and the matrix $M$ then can
simply be chosen as having the vectors $v^W_i$ as its columns, i.e. $M = (v^W_{ij})$.
By associating to
$\mathcal{V}$ the structure of some appropriate fine-graded module, the matrix $M$
becomes a monomial matrix for which syzygies can be computed.

\begin{definition}
\begin{enumerate}[(i)]
\item A {\em reflexive model} for $\mathcal{V}$ is an inclusion $\mathcal{V}
\hookrightarrow (\Z^r, \leq)$ for some $r > 0$ such that that its image in $\Z^r$ is
an $\lcm$-lattice.
\item A set of {\em generating flags} of $\mathcal{V}$ is a set of tuples
$\{E^1_1 \subsetneq \cdots \subsetneq E^1_{n_1}\}, \dots, \{E^r_1 \subsetneq \dots
\subsetneq E^r_{n_r}\} \subset \mathcal{V}$ such that $E^i_{n_i} = \mathbf{V}$ for
every $i$ and $\mathcal{V}$ is the set of all intersections among the $E^i_j$.
\end{enumerate}
\end{definition}

Let $E^i_j$ be any set of generating flags, then we associate to each of the flags
a tuple of integers $\underline{k}^i := (k^i_1 < \dots < k^i_{n_i})$.
This data defines a reflexive model, where we map every $W \in \mathcal{V}$ to the
tuple $\underline{k}_V := (\min\{k^i_j \mid W \subset E^i_{k^i_j}\} \mid i = 1,
\dots, r)$. As easily
can be seen, this reflexive model gives rise to a reflexive fine-graded module $E$
over the polynomial ring $S = k[x_1, \dots, x_r]$ which is given by filtrations
\begin{equation*}
E^i(j) =
\begin{cases}
0 & \text{ if } j < k^i_1, \\
E^i_l & \text{ if } k^i_l \leq j < k^i_{l + 1}, l < n_i, \\
\mathbf{V} & \text{ if }  n_i \leq j.
\end{cases}
\end{equation*}
$E$ has an embedding into the free module $S(-\underline{k}_{\min})^{\dim
\mathbf{V}}$, where
$\underline{k}_{\min} = (k^1_1, \dots, k^r_1)$, and the module $F_0$ is given by
the direct sum $\bigoplus_{W \in \mathcal{V}} S(-\underline{k}_V)^{f_V}$, where
$f_V$ is the free dimension of $V$. So, we have
\begin{equation*}
0 \longrightarrow K_0 \longrightarrow \bigoplus_{W \in \mathcal{V}}
S(-\underline{k}_V)^{f_V} \overset{\bar{M}}{\longrightarrow}
S(-\underline{k}_{\min})^{\dim \mathbf{V}},
\end{equation*}
where $\bar{M}$ is a monomial matrix whose entries are of the form
$(v^W_{ij} x^{\underline{k}_V - \underline{k}_{\min}})$, where the $v^W_{ij}$ are the
corresponding entries of the matrix $M$. The image of $M$ then is the module $E$.
So the only effect seen by the choice of the reflexive model for $\mathcal{V}$ are
the number of variables in the ring $S$ and the degrees of the monomials in
$\bar{M}$ and the subsequent matrices in the resolution, whereas the coefficients in
$\bar{M}$ are precisely the entries of the matrix $M$.

\addtocontents{toc}{\vskip4mm}
\addtocontents{toc}{\bf References \hfill \thepage}

\end{document}